\documentclass[twoside]{article}
\usepackage[utf8]{inputenc}
\usepackage{amsmath}
\usepackage{amsfonts}
\usepackage{amssymb}
\usepackage{graphicx}
\usepackage{color}
\usepackage{float}
\usepackage[left=2cm,right=2cm,top=2cm,bottom=2cm]{geometry}
\usepackage{fdsymbol}
\newenvironment{Proof}{\noindent{\sc Proof.}}{\qed}
\newtheorem{theorem}{Theorem}
\newtheorem{lemma}{Lemma}[section]

\newcommand{\qed}{\hfill$\Box$\par\medskip}

\def\bhag#1{\noindent
\setcounter{equation}{0}
\section{#1}
}

\def\RR{{\mathbb R}}
\def\CC{{\mathbb C}}
\def\ZZ{{\mathbb Z}}

\def\PPI{{{\rm I}\kern-1pt\Pi}}

\def\b #1;{{\bf #1}}
\def\x{{\bf x}}

\def\z{{\bf z}}

\def\F{{\cal F}}

\def\M{{\mathcal M}}

\def\E{{\cal E}}

\def\be{\begin{equation}}
\def\ee{\end{equation}}
\def\bea{\begin{eqnarray}}
\def\eea{\end{eqnarray}}

\def\donchitre#1#2{\vskip 6.5cm\noindent
\parbox[t]{1in}{\special{eps:#1.eps x=6.5cm y=5.5cm}}
\hbox to 7cm{}\parbox[t]{0.0cm}{\special{eps:#2.eps x=6.5cm y=5.5cm}}}

\def\CC{{\mathbb C}}

\begin{document}
\title{Order of uniform approximation by polynomial interpolation in the complex plane and beyond}
\author{Charles K. Chui$^{}$\thanks{{Charles lives in Menlo Park, CA and is affiliated with the Statistics Department of Stanford University, CA. His research is partially supported by the U.S. Army Research Office under ARO Grant \#W911NF2110218.} \textsf{Email:} ckchui@stanford.edu}
\ and
Lefan Zhong$^{}$\thanks{{Independent Researcher.} \textsf{Email:} lefan$\_{\rm z}$@yahoo.com}}

\date{}

\maketitle
\begin{abstract}
Let $X$ denote a compact set in the complex plane $\CC$ and ${\bf {z}}_n:=\{z_{n,j}\}^{n}_{j=0}$ be a family of points that lie on $X$. Then the norm $\|A_{{\bf {z}}_n}\|$ of the operator $A_{{\bf {z}}_n}$ that maps the Banach space $C(X)$ to the space $\Pi_n$ of polynomials of degree $\le n$, defined by $A_{{\bf {z}}_n}f = p_n$, where $p_n\in \Pi_n$ is the polynomial that interpolates $f\in C(X)$ at ${\bf {z}}_n$, is called the Lebesgue constant of the family ${\bf {z}}_n$. While the bulk of this paper is on the study of Lagrange polynomial interpolation at the Fej\'er points ${\bf z}_n^{*}:=\{z_{n,j}^{*}\}^{n}_{j=0}$, and their proper adjustment, that lie on an open arc $X= \gamma$ in $\CC$ which does not cross itself, the general spirit of our presentation carries over to other settings of polynomial representations of functions on a simple closed piece-wise smooth curve $X=\Gamma$, which is the boundary of a Jordan domain $D\subset \CC$; and instead of polynomial interpolation, $A_{{\bf {z}}_n}$ is replaced by the linear operator $A_{F_n}$, that maps the (generalized) Hardy space $H^{p}(D)$, for $p\ge 1$, of analytic functions $f$ in $D$ with non-tangential limit $f^{*} \in L_p(\Gamma)$, to the $n$-th partial sum $S_n(\cdot; f)$ of the Faber series representation of $f$ in $D$. The importance of the Lebesgue constant $\|A_{F_n}\|$ is that if $P_n^{*}\in \Pi_n$ is the best polynomial approximant of $f$ in $H^{p}(D)$, with approximation error $\epsilon_n (f; \Gamma):= \| f - P_n^{*}(\cdot; f)\|_{H^{p}}$, then $S_n(\cdot;f)$ can be used to replace $P_n^{*}$ with error of approximation $\|f-S_n(\cdot;f)\|_{H^{p}}\le (1+\|A_{F_n}\|)\,\epsilon_n (f; \Gamma)$. This consideration, along with the discussion of related problems and results inspired by the work of Prof. J. Korevaar, will be presented in the final section of the present paper. 

For Lagrange polynomial interpolation on open arcs $X=\gamma$ in $\CC$, it is well-known that the Lebesgue constant for the family of Chebyshev points ${\bf{x}}_n:=\{x_{n,j}\}^{n}_{j=0}$ on $[-1,1]\subset \RR$ has growth order of $O(\log(n))$. The same growth order was shown in \cite{ZZ} for the Lebesgue constant of the family ${\bf {z^{**}_n}}:=\{z_{n,j}^{**}\}^{n}_{j=0}$ of some properly adjusted Fej\'er points on a rectifiable smooth open arc $\gamma\subset \CC$. On the other hand, in our recent work \cite{CZ2021}, it was observed that if the smooth open arc $\gamma$ is replaced by an $L$-shape arc $\gamma_0 \subset \CC$ consisting of two line segments, numerical experiments suggest that the Marcinkiewicz-Zygmund inequalities are no longer valid for the family of Fej\'er points ${\bf z}_n^{*}:=\{z_{n,j}^{*}\}^{n}_{j=0}$ on $\gamma$, and that the rate of growth for the corresponding Lebesgue constant $L_{{\bf {z}}^{*}_n}$ is as fast as $c\,\log^2(n)$ for some constant $c>0$.

The main objective of the present paper is 3-fold: firstly, it will be shown that for the special case of the $L$-shape arc $\gamma_0$ consisting of two line segments of the same length that meet at the angle of $\pi/2$, the growth rate of the Lebesgue constant $L_{{\bf {z}}_n^{*}}$ is at least as fast as $O(\log^2(n))$, with $\lim\sup \frac{L_{{\bf {z}}_n^{*}}}{\log^2(n)} = \infty$; secondly, the corresponding (modified) Marcinkiewicz-Zygmund inequalities fail to hold; and thirdly, a proper adjustment ${\bf z}_n^{**}:=\{z_{n,j}^{**}\}^{n}_{j=0}$ of the Fej\'er points on $\gamma$ will be described to assure the growth rate of $L_{{\bf z}_n^{**}}$ to be exactly $O(\log^2(n))$.
\end{abstract}
 \vskip 10pt
 \begin{center} \bf{Dedicated to Prof. Jaap Korevaar on the occasion of his 100-th birthday!}
 \end{center}
\bhag{Introduction and results}
Representation of functions by polynomials has been one of the most well studied research areas in classical analysis. If a complicated function $f(x)$ is adequately represented by some polynomial  $p_n(x)$ of degree $n$, then the tasks of manipulating, understanding and editing the information contained in $f(x)$ are possible and easily accomplished by processing the polynomial $p_n(x)$ instead. Such tasks include differentiation, integration, visualization and geometric design based on the function $f(x)$. An effective measurement of the adequacy in polynomial representation is the sup-norm approximation, denoted by $ \|f-p_n\|_{\infty}$. For instance, let $\Pi_n$ denote the space of polynomials of degree $\le n$ and consider real-valued bounded measurable functions $f(x)$ on the interval $[-1, 1]\subset \RR$, then $P_n(\cdot; f) \in \Pi_n$ is the best approximation of $f$, if 
\begin{equation} \label{best approx}
\| f - P_n(\cdot; f)\|_{\infty} = E_n(f; [-1,1]) := \inf \{\|f - p_n\|_{\infty}: p_n\in \Pi_n\}.
\end{equation}
While $f$ could be a very complicated function, it may still be a function in $C^{m}[-1, 1]$ for some integer $m\in \ZZ_{+}:=\{0,1, \cdots\}$. For such functions, recall the well-known theorem of Jackson (see \cite{J}) that assures the fast rate of uniform convergence of $P_n(\cdot; f)$ to $f$ on $[-1, 1]$, namely:
\begin{equation} \label{approx order}
E_n(f; [-1,1])= O(\frac{1}{n^{m}})\omega(\frac{1}{n}; f^{(m)}),
\end{equation}
where $\omega(\frac{1}{n}; f^{(m)})$ denotes the uniform modulus of continuity of $f^{(m)}$. Thus, it is desirable to replace $f$ by $P_n$ for carrying out the tasks as mentioned above. However, it is not an easy task to derive the polynomial $P_n$ that satisfies (\ref{best approx}). On the other hand, it is relatively easy to compute the polynomial $p_n = p_n(\cdot;f,{\bf x}_n)$ that interpolates the set $\{f(x_{n,j})\}$ of discrete samples of $f$ at any desired sampling locations (commonly called interpolation nodes) ${\bf x}_n := \{x_{n,j}\}$, in that $p_n(x_{n,j}) = f(x_{n,j})$ for $j=0,\cdots,n$. 

Here and throughout this paper, for any family of points ${\bf {z}}_n:=\{z_{n,j}\}^{n}_{j=0}$ that lie on any compact set $X\subset \CC$, including the interval $[-1,1]$, we always assume that for each $n\in \ZZ_{+}$, the points $z_{n,0},\cdots,z_{n,n}$ are distinct, and for this family ${\bf {z}}_n$ of points on $X$, we also introduce the notion of ``Lagrange fundamental polynomials" defined by: 
\begin{equation} \label{lagrange poly}
\ell_{n,j}(z) :=\prod_{k\neq j, \,0\le k \le n} \dfrac{z-z_{n,k}}{z_{n,j} - z_{n,k}},\, \,\,\,\,\,\, j=0,\cdots,n.
\end{equation}

Returning to the computation of the polynomials $p_n = p_n(\cdot;f,{\bf x}_n)$ that interpolate $f$ at $\{x_{n,j}\}$, we may simply pre-compute the data-independent ``Lagrange fundamental polynomials" $\{\ell_{n,j}(x)\}_{n,j} \subset \Pi_n$, and use the samples $f(x_{n,j})$ of the given function $f(x)$ to obtain
\begin{equation} \label{Lag int}
p_n(x;f,{\bf x}_n)=\, \sum^{n}_{j=0} {f(x_{n,j})\, \ell_ {n,j}(x)}.
\end{equation}
The challenge of this approach of replacing the best approximation polynomial by polynomial interpolation is the selection of interpolation nodes ${\bf x}_n =\{x_{n,j}\}$. Indeed, it is well-known that the choice of equally spaced sampling nodes ${\bf x}_n^{r} :=\{ x^{r}_{n,j} = -1 +\frac{2j}{n}\}^n_{j=0}$ is not recommended due to the so-called Runge's phenomenon \cite{RUN1901}, meaning that even for functions in $C^{\infty}[-1, 1]$, such as the analytic function $g(x):= \frac{1}{1+25x^2}$, the error of uniform approximation $\|g - p_n(\cdot;g,{\bf x}^{r}_n)\|_{\infty}$ of $g$ by the interpolating polynomial increases (instead of being decreasing) with increasing values of the polynomial degrees $n\in \ZZ_{+}$ (see \cite{B} pages 22 -- 26). In view of (\ref{best approx}) and (\ref{approx order}), it is certainly unacceptable to replace the polynomial $P_n(\cdot; f)$ by $p_n(x;f,{\bf x}^{r}_n)$ to represent any continuous function on $[-1, 1]$. 

So what could be a good selection criterion of interpolation nodes? Let us address this question for the more general setting of interpolation nodes ${\bf {z}}_n:=\{z_{n,j}\}^{n}_{j=0}$ that lie on a compact set $X$ in the complex plane $\CC$, including the interval $[-1, 1]$ as discussed above. For a function $f\in C(X)$, suppose that the error of uniform approximation
\begin{equation} \label{poly approx}
\epsilon_n (f; X):= \| f - P_n(\cdot; f)\|_{\infty}  
\end{equation}
of $f$ by some polynomial $P_n(\cdot;f)\in \Pi_n$ is $O(\frac{1}{n^{\beta}})$ for some $\beta>0$ which could be arbitrarily small. Is it worthwhile to replace the polynomial $P_n(\cdot; f)$, which is not easy to derive, by the polynomial $p_n = p_n(\cdot;f,{\bf z}_n)$ that interpolates the discrete data samples $\{f(z_{n,j})\}^{n}_{j=0}$? To answer this question, let $w_{n,j} := P_n(z_{n,j}; f)$ and consider the polynomial $p_n(\cdot;P_n,{\bf z}_n)$ that interpolates the discrete data $\{w_{n,j}\}^{n}_{j=0}$. Then since the polynomial $q_n(z) := p_n(z;P_n,{\bf z}_n) - P_n(z;f)$ is a polynomial of degree $\le n$ and has $n+1$ distinct zeros $z_{n,0},\cdots,z_{n,n}$, it must be the zero polynomial; that is, the two polynomials $p_n(z;P_n,{\bf z}_n)$ and $ P_n(z;f)$ are identical. Now, for each $z\in X$, by applying (\ref{lagrange poly}), we have
\begin{eqnarray}\label{Lag error}
|f(z)- p_n(z;f,{\bf z}_n)| &\le& |f(z)- P_n(z; f)| + |P_n(z;f)-p_n(z;f,{\bf z}_n)| \nonumber\\ 
&=& |f(z)- P_n(z; f)| + |p_n(z;P_n,{\bf z}_n) - p_n(z;f,{\bf z}_n)| \nonumber\\
&=&|f(z)- P_n(z; f)| + |\sum^{n}_{j=0} (w_{n,j} - f(z_{n,j})\, \ell_ {n,j}(z)| \nonumber\\
&\le& |f(z)- P_n(z; f)|  + \sum^{n}_{j=0} |P_n(z_{n,j}; f) - f(z_{n,j})|\, |\ell_ {n,j}(z)| \nonumber\\
&\le& \epsilon_n(f; X)(1 + \max_{z\in X} \, \sum^{n}_{j=0} |\ell_ {n,j}(z)|), 
\end{eqnarray}
where $\Lambda(z):=\sum^{n}_{j=0} |\ell_ {n,j}(z)|$ is called the Lebesgue function in the literature (see, for example $\cite{R1969}$), and $\max_{z\in X} \, \Lambda (z)$ is often called the ``Lebesgue constant" (see $\cite{R1969}$);
but for the reason to be discussed in the Final Remarks of this paper, we prefer to define the Lebesgue constant $L_{{\bf {z}}_n}$ of ${\bf {z}}_n$ in terms of the norm of the bounded linear operator that maps the Banach space $C(X)$ to $\Pi_n$ by Lagrange polynomial interpolation at ${\bf {z}}_n$ in this paper. As a consequence of (\ref{Lag error}), the ``selection criterion of interpolation the nodes ${\bf {z}}_n:=\{z_{n,j}\}^{n}_{j=0}$" should at least be $\epsilon_n (f; X) L_{{\bf {z}}_n} \rightarrow 0$ for $n\rightarrow \infty$. 

Returning to Jackson's theorem on the order of best polynomial approximation $E_n(f; [-1,1])= O(\frac{1}{n^{m}})\omega(\frac{1}{n}; f^{(m)})$ for functions $f\in C^{m}[-1,1]$, we must ensure that the choice of interpolation nodes  ${\bf x}_n:=\{x_{n,j}\}^{n}_{j=0}$ must at least satisfy $\frac{L_{{\bf {z}}_n}}{n^m}\omega(\frac{1}{n}; f^{(m)})\rightarrow 0$ for increasing values of $n$. While the family  ${\bf x}_n^{r} :=\{ x^{r}_{n,j} = -1 +\frac{2j}{n}\}^n_{j=0}$ of equally spaced nodes fails to satisfy this selection criterion, the family ${\bf x}_n^{c} :=\{ x^{c}_{n,j}\}^{n}_{j=0}$, with $x^{c}_{n,j}=\cos\big({{2j + 1}\over{2(n+1)}}\pi\big),$ for $j=0,1,\cdots,n$, called the ``Chebyshev points" of the interval $[1, -1]$, can be used as interpolation nodes for any function $f\in Lip_{\alpha}[-1,1]$ for any arbitrarily small $\alpha > 0$, for which $\omega(\frac{1}{n};f) = O(\frac{1}{n^{\alpha}})$. The reason is that the Lebesgue constant  $L_{{\bf x}^{c}_n}$ of the family ${\bf x}_n^{c}$ satisfies $L_{{\bf x}^{c}_n}= O(\log(n))$, and that $\frac{\log(n)}{n^{\alpha}}\rightarrow 0.$ For Chebyshev points that lie on the interval $[1, -1]$, there is a vast amount of literature on the upper and lower bounds, and even on asymptotic rate of growth, of $L_{{\bf x}^{c}_n}= O(\log(n))$. For the interested reader, we list the work of Faber $\cite{F1914}$, Bernstein $\cite{BER1931}$, Erd\"os $\cite{ERD1958}$ and V\'ertesi $\cite{VER1990}$ in the chronological order, and the more recent review paper $\cite{Ib2016}$. As to the order of uniform approximation of continuous functions on a rectifiable curve in the complex plane $\CC$, the interested reader is referred to $\cite{D1977}$, $\cite{Belyi1977}$ and  $\cite{Gaier1987}$.

To consider the Lebesgue constant of a general family ${\bf {z}}_n:=\{z_{n,j}\}^{n}_{j=0}$ of interpolation nodes that lie on a compact set $X\in \CC$, we always assume that for each $n$, the $n+1$ points of the set $\{z_{n,0},\cdots,z_{n,n}\}$ are distinct. Consider the Banach space $C(X)$ with the sup-norm over $X$ and the bounded linear operator $A_{{\bf {z}}_n}$ that maps  $C(X)$ to $\Pi_n$, defined by 
\begin{equation} \label{operator}
(A_{{\bf {z}}_n} f)(z) := p_n(z;f, {\bf {z}}_n), 
\end{equation}
for each $f\in C(X)$, with operator norm 
\begin{equation} \label{operator norm}
\|A_{{\bf {z}}_n}\|:= \it{\sup}\,\{\|p_n(.\,;f)\|_{\infty}: f\in C(X),\,\, \|f\|_{\infty} = 1 \}.
\end{equation}
Then the Lebesgue constant $L_{{\bf {z}}_n}$ is defined by
\begin{equation} \label{Lc def}
L_{{\bf {z}}_n}:= \|A_{{\bf {z}}_n}\|.
\end{equation}
For the convenience of the reader, we will derive the following well-known definition of $L_{{\bf {z}}_n}$ in terms of the Lebesgue function and give another well-known formula for computation.  
\begin{lemma} \label{leb const}
Let ${\bf {z}}_n:=\{z_{n,j}\}^{n}_{j=0}$ be a family of nodes in $X$. Then the Lebesgue constant $L_{{\bf {z}}_n}$ defined in (\ref{Lc def}) can be formulated as:
\begin{equation}\label{Lcformulation}
L_{{\bf {z}}_n}= \max_{z\in X} \, \sum^{n}_{j=0} |\ell_ {n,j}(z)|,
\end{equation}
and can be computed by applying the formula:
\begin{equation} \label{ineq:Lebesg1}
L_{{\bf {z}}_n} = \max_{z\in X} \sum_{k=0}^n {\frac{|\omega_n (z)|} {{|\omega_n^\prime (z_{n,k}) } (z-z_{n,k} ) |}}\,\,,
\end{equation}
where 
\begin{equation} \label{Omega0}
\omega_n(z) := \prod_{0\le k \le n} (z - z_{n,k}).
\end{equation}
\end{lemma}
\begin{Proof}
To prove (\ref{Lcformulation}), let $z \in X$ be fixed and consider any function $f^*\in C(X)$ that satisfies $\|f^*\|_{\infty} = 1$ and 
$$
f^*(z_{n,j}) = { \frac{\overline{\ell_{n,j}(z)}}{\ell_{n,j}(z)}},
$$
so that it follows from (\ref{Lag int}) that
$$\|p_n(.;f^{*})\|_{\infty}=\,\max_{z\in X} \sum^{n}_{j=0} |\ell_ {n,j}(z)|,$$
and by the definition of $\|A_n\|$ in (\ref{operator norm}), we have 
\begin{equation} \label{lower bd}
\|A_{{\bf {z}}_n}\|\ge \|p_n(.\,;f^{*})\|_{\infty} = \max_{z\in X} \sum^{n}_{j=0} |\ell_ {n,j}(z)|.
\end{equation}
On the other hand, from the definition of $A_{{\bf {z}}_n}$ in (\ref{operator}) and the choice of $f^{*}$, it is clear that
\begin{equation} \label{upper bd}
\|A_{{\bf {z}}_n}\|\le \max_{z\in X} \, \sum^{n}_{j=0} |\ell_ {n,j}(z)|.
\end{equation}
Therefore, by combining (\ref{lower bd}) and (\ref{upper bd}), it follows from the definition (\ref{Lc def}) of the Lebesgue constant that 
\begin{equation} \label{operator norm 2}
L_{{\bf {z}}_n}=\|A_{{\bf {z}}_n}\|\ = \max_{z\in X} \, \sum^{n}_{j=0} |\ell_ {n,j}(z)|.
\end{equation}
To prove (\ref{ineq:Lebesg1}), we re-write (\ref{Omega0}) as 
$$\omega_n(z) = (z-z_{n,k})\, \prod_{j\neq k, \,0\le j \le n}(z - z_{n,j})\, ,$$
for each fixed $k$, take the derivative of this product (of two functions), and evaluate at $z = z_{n,k}$, to yield 
$$\omega_n^\prime (z_{n,k}) =  \prod_{j\neq k, \,0\le j \le n}(z_{n,k} - z_{n,j})\,.$$
This formula enables us to re-write the  Lagrange fundamental polynomials $\ell_{n,k}(z)$ in (\ref{lagrange poly}) as 
$$\ell_{n,k}(z)= {\frac{\omega_n (z)} {\omega_n^\prime (z_{n,k})  (z-z_{n,k} )}}\,.$$
This completes the proof of the lemma. 
\end{Proof}
\vskip .1 in

To extend our discussion of the Chebyshev points ${\bf x}_n^{c} :=\{ x^{c}_{n,k}\}^{n}_{k=0}$ on the interval $[-1, 1]$, namely: $x^{c}_{n,k}=\cos\big({{2k + 1}\over{2(n+1)}}\pi\big)$, where $k=0,\dots,n$, to the study of interpolation nodes ${\bf {z}}_n:=\{z_{n,k}\}^{n}_{k=0}$ that lie on a rectifiable open arc $\gamma \subset \CC$ which does not cross itself, let $\psi_\gamma(w)$ denote the conformal map of $|w| > 1$ to $\CC^{*}\backslash \gamma$, with $\psi_\gamma(\infty) = \infty$ and $\psi_\gamma^{'}(\infty) > 0$, where $\CC^{*}$ denotes the extended complex plane, and let $\psi_\gamma(w)$ be extended continuously to the unit circle, so that $\psi_\gamma(w)$ maps  $|w|=1$ onto $\gamma$. Then for any fixed real number $\theta_n$, the family ${\bf z}^{*}_n := {\bf z}^{*}_{n,\theta} := \{z^{*}_{n,k}\} \subset \gamma$, where 
\begin{equation} \label{Fej points}
z^{*}_{n,k} = \psi_\gamma(e^{i{\frac{2 k\pi}{n + 1} + \theta_n}}),\hskip 10pt k = 0,1,\cdots,n,
\end{equation}
defined in \cite{FE1918} is called the family of Fej\'er points on $\gamma$ with parameter $\theta_n$. For example, for $\gamma = [-1, 1] \subset \RR$, the conformal map $\psi_X := \psi_{[-1,1]}$ of $|w| > 1$ to $\CC \backslash [-1, 1]$ is given by $\psi_X(w) = \frac{1}{2}(w + \frac{1}{w})$ with continuous extension to $|w|= 1$ being $\psi_X(e^{i\theta}) = \cos\,\theta$. Hence, by selecting $\theta_n =\frac{\pi}{2(n + 1)}$, the Fej\'er points $\{\psi_X\big(e^{i(\frac{(2k+1)\pi}{2(n + 1)})}\big)\}^{n}_{k=0}$ are precisely the Chebyshev points $\{\cos\big({{(2k + 1})\pi\over{2(n+1)}}\big)\}^{n}_{k=0}$ on $[-1, 1]$. Hence, it seems to be quite natural to consider Fej\'er points on $\gamma\subset \CC$ as an extension of Chebyshev points on $[-1, 1]\subset\RR$. More on this topic will be discussed in the final section.\\

Before going into the discussion of the Lebesgue constants of polynomial interpolation at the Fej\'er points $\{z^{*}_{n,k}\}^{n}_{k=0}$ that lie on an open rectifiable arc $\gamma\subset \CC$ which does not cross itself, let us recall the paper \cite{CUR1935} by J. H. Curtiss, in which it was shown in Theorem III b that for a bounded Jordan domain $D$ with smooth boundary curve $\Gamma$, the order of uniform approximation by Lagrange interpolation polynomials $p_n(z;f,{\bf z}^{*}_n)\in \Pi_n$, where ${\bf z}^{*}_n$ denotes the family of Fej\'er points that lie on $\Gamma$, is $\log(n)$ multiple of the best order of uniform approximation from $\Pi_n$. Hence, in view of the argument in (\ref{poly approx}) and (\ref{Lag error}), Curtiss' result implies that the order of growth  of the Lebesgue constant $L_{{\bf {z}}^{*}_n}$ of Fej\'er points on a smooth Jordan curve is $O(\log(n))$. This result of $O(\log(n))$ growth order was later shown in \cite{ZZ} to be valid for the Lebesgue constant $L_{{\bf {z}}^{**}_n}$ of properly adjusted Fej\'er points $\{z^{**}_{n,j}\}^{n}_{j=0}$ that lie on an open smooth arc $\gamma\subset \CC$ that does not cross itself. The obvious question to ask is if this growth order holds for piece-wise smooth open arcs. The main objective of the present paper is to address this question.

To convey the motivation of our research on this problem, let us first digress our discussion to mentioning the well-known Marcinkiewicz-Zygmund inequalities $\cite{MZ}$ that are instrumental to the study of trigonometric polynomial interpolation, as follows: For $1<p< \infty$, there exist constants $0< c_1 \le c_2$, such that
\begin{equation}\label{MZ ineq}
{{c_1}\over{n+1}} \sum^n_{k=0} \Big|P_n \Big(e^{i{{2\pi k}\over{{n+1}}}}\Big)\Big|^p \le {1\over{2\pi}}\int_0^{2\pi} \Big|P_n (e^{i\theta}) \Big|^p  d\theta \le {{c_2}\over {n+1}} \sum^n_{k=0} \Big|P_n \Big(e^{i{{2\pi k}\over{{n+1}}}}\Big)\Big|^p,
\end{equation}
for all polynomials $P_n(z)\in \Pi_n$ and for all $n\in \ZZ_{+}$, where $z=e^{i\theta}$. To extend the trigonometric polynomials $P(e^{i\theta})$ to algebraic polynomials $P_n(z)$ for $z\in \CC$, we replace the points $z_{n,k}=e^{i\frac{2\pi k}{n+1}}$ by a more general family of interpolation notes ${\bf {z}}_n:=\{z_{n,k}\}^{n}_{k=0}$ that lie on the unit circle $|z|=1$, and consider the extension of (\ref{MZ ineq}) to 
\begin{equation}\label{gMZ ineq}
{{c_1}\over {n+1}}\sum_{k=0}^n \big|P_n (z_{n,k}) \big|^p \le {1\over 2\pi}\int_{0}^{2\pi}\big|P_n(e^{i\theta})\big|^p  d\theta \le {{c_2}\over {n+1}} \sum_{k=0 }^n \big|P_n (z_{n,k} ) \big|^p.
\end{equation}
Of course the constants $c_1, c_2$ in (\ref{gMZ ineq}) are different from those in (\ref{MZ ineq}) and must be independent of ${\bf z}_n$. The validity of this extension was shown in our previous work \cite{CZ1999} under two assumptions: firstly, the separation condition
\begin{equation}\label{sep cond}
 \min_{0\le k<n} |z_{n,k+1} - z_{n,k} | \ge {{c_3}\over{n+1}},
\end{equation}
for some arbitrarily small constant $c_3>0$; and secondly, the $A_{p}$-weight condition of $\omega_n(z) = \prod_{0\le j \le n}(z - z_{n,j})$ defined in (\ref{Omega0}), namely:
 \begin{equation}\label{A_p 1}
\sup_{t_0<t_1 } \Big \{ {1\over{t_1-t_0}} \int_{t_0}^{t_1} | \omega_n ((1+{{1}\over {n+1}}) e^{it} ) |^p dt \Big \}^{1\over p} \Big \{{1\over {t_1-t_0}} \int_{t_0}^{t_1} |\omega_n ((1+{{1}\over {n+1}}) e^{it}) |^{-q} dt \Big\}^{1\over q}\le c_4,
\end{equation}
for some constant $c_4>0$, where $p, q > 0$ and $\frac{1}{p} + \frac{1}{q} = 1$. Indeed, it was shown in \cite{CZ1999} that the sufficient conditions (\ref{sep cond}) and (\ref{A_p 1}) are also necessary for the validity of the extended Marcinkiewicz-Zygmund inequalities in (\ref{gMZ ineq}) for the unit circle $|z| =1$.

As to the consideration of Marcinkiewicz-Zygmund inequalities for interpolation nodes ${\bf x}_n =\{z_{n,k}:k=0,1,\dots,n\}$ that lie on an open arc $\gamma\subset \CC$, we proposed in our recent paper \cite{CZ2021} the following ``modified" version that depends on the conformal map $\psi_\gamma$ from $|w| > 1$ onto $\CC^{*}\backslash \gamma$ as discussed above. Let $\Gamma_n$ be the level curves $\{\psi(w): |w|= 1 +\frac{1}{n+1}\}$
and denote by $\text{dist}(z, \Gamma_n)$ the distance from $z\in \gamma$ to the level curve $\Gamma_n$. Then for the modified version, we proposed to consider $\text{dist}(z_{n,k}, \Gamma_n)$ for $k=0,\cdots,n$, instead of the distances between the neighboring points of the $\{z_{n,k}:k=0,1,\dots,n\}$, as follows.\\

For $1<p< \infty$, there exist constants $0< c_1 \le c_2$, such that
\begin{equation} \label{mod MZ ineq}
c_1 \sum^n_{k=0} |P_n (z_{n,k} ) |^p \text{dist}(z_{n,k},\Gamma_n) \le \int_\gamma |P_n (z) |^p |dz| \le  c_2 \sum^n_{k=0} |P_n (z_{n,k} ) |^p \text{dist}(z_{n,k},\Gamma_n ),
\end{equation}
for all polynomials $P_n(z)\in \Pi_n$ and for all $n\in \ZZ_{+}$. For the special case of $\gamma$ consisting of two 2 line segments in $\CC$, it was shown in \cite{CZ2021} that validity of the modified Marcinkiewicz-Zygmund inequalities in (\ref{mod MZ ineq}) is equivalent to the totality of the $A_{p}$-weight condition of $\omega_n(z) = \prod_{0\le j \le n}(z - z_{n,j})$ defined in (\ref{A_p 1}) and the $c_0$-separation condition, defined by 
\begin{equation} \label{ineq:sep_c0}
\min_{0\le k < j\le n}   {{|z_{n,j}-z_{n,k}|} \over {\min \big(\text{dist}(z_{n,j},\Gamma_n ),\text{dist}(z_{n,k},\Gamma_n ) \big)}} \ge c_0>0.
\end{equation}

Returning to the study of the growth order of Lebesgue constants for a rectifiable smooth open arc $\gamma$ that does not cross itself,  we will show that the $O(\log(n))$ is no longer valid for piece-wise smooth open arcs (with at least one corner), even for adjusted Fej\'er points. In this regard, since such claims are negative results, it is sufficient to consider the special case of $\gamma_0$ that consists of two line segments of equal length, meeting at the origin of $\CC$ at an angle of $\frac{\pi}{2}$ and symmetric with respect to the real axis, namely:
\begin{equation} \label{definegamma}
\gamma_{0} := \{(-1+i)x: \,x\in [0,27^\frac{1}{4}/\sqrt{2}]\}\cup\{(-1-i)x:\,x\in [0,27^{1\over 4}/\sqrt{2}]\}.
\end{equation}
The conformal map of $|w| > 1$ onto $\CC\backslash\gamma_0$, is given by
\begin{equation} \label{ineq:Conform0}
\psi_{0} (w) =\Big(w-{\frac{1} {w}}\Big) \Big({\frac{w-1}{w+1}}\Big)^{\frac{1}{2}},
\end{equation}
and extended continuously to the unit circle, so that $\psi_0(w)$ maps $|w|=1$ onto $\gamma_{0}$.\\

To select the most suitable Fej\'er points $z^{*}_{n,k},\,\, k=0,\cdots,n$,\, on $\gamma_0$, it is necessary to consider suitable choices of the proper ``rotation" of the equally spaced points on the unit circle $|w|=1$ in order to attain sufficiently uniform spacing of the points $z^{*}_{n,k}$ on $\gamma_0$, and also to avoid ``collapses" among the points $z^{*}_{n,k}$. Hence, for the sake of ``symmetry", we consider the following selection.\\

For every even positive integer $n$, we choose the $n+1$ roots of $w^{n+1} - 1 = 0$ as $\{w_{n,k}\}_{k = 0}^n$, while for every odd positive integer $n$, the $n+1$ roots of $w^{n+1} + 1 = 0$ are chosen to be the points on the unit circle. The reason for this choice is that since the open arc $\gamma_0$ is symmetric with respect to the x-axis, the choice of $\{w_{n,k}\}_{k = 0}^n$ should also be symmetric respect to the x-axis, so that in view of the normalization condition $\psi_0(\infty)>0$ of the exterior conformal map $\psi_0$, the Fej\'er points $\{\psi_0(w_{n,k})\}_{k = 0}^n \subset \gamma_0 $ have the same symmetry property. More specifically,

\vskip .1in
\noindent(1) for $ 0 \le k \le \lfloor{n/2}\rfloor$, we consider
\begin{equation}\label{ineq:theta0}
\theta_{n,k} = \left\{
\begin{array}{ll}
      {{2k\pi}\over{n+1}}, & n = 2m,  \\
      \\
      {{(2k+1)\pi}\over{n+1}}, & n = 2m+1; \\
\end{array} 
\right. 
\end{equation}
\noindent(2) for $\lfloor{n/2}\rfloor < k \le n$, we consider
\begin{equation}\label{ineq:theta11}
\theta_{n,k} = -\theta_{n,2\lfloor{n/2}\rfloor + 1 - k}, \hskip 10pt \lfloor{n/2}\rfloor < k \le n,\\
\end{equation}
and set $w_{n,k} := e^{i\theta_{n,k}}$, \,for $k=0,\cdots, n$. In the following two theorems of this paper, we only consider the  Fej\'er points ${\bf z}^{*}_n :=\{z^*_{n,k}\}_{k=0}^n \subset \gamma_0$, defined by
\begin{equation}\label{ineq:Fejer0}
z^*_{n,k} := \psi_0(e^{i\theta_{n,k}}), \hskip 15pt k = 0,\cdots,n.
\end{equation}

\begin{theorem} \label{theorem:JustFejer}
Let $\gamma_0$ be the $L$-shape arc in $\CC$ defined in $(\ref{definegamma})$ and  ${\bf z}^{*}_n =\{z^*_{n,k}\}_{k=0}^n$ be the family of Fej\'er points on $\gamma_0$ defined in (\ref{ineq:Fejer0}). Then there exists a constant $c_1 > 0$, such that 
\begin{equation} \label{ineq:Lebesg0}
c_1 \log^2(n) \le L_{{\bf {z^*_n}}},
\end{equation}
for all $n\in \ZZ_{+}$. On the other hand, 
\begin{equation} \label{ineq:LebesgInfty}
\lim\sup_{n\rightarrow\infty} \frac{L_{{\bf {z^*_n}}}}{\log^2(n)} = \infty.
\end{equation}

\end{theorem}

\begin{theorem} \label{theorem:MZ_fail}
Let $\gamma_0$ be the $L$-shape arc in $\CC$ defined in $(\ref{definegamma})$ and ${\bf z}^{*}_n =\{z^*_{n,k}\}_{k=0}^n$ be the family of Fej\'er points on $\gamma_0$ defined in (\ref{ineq:Fejer0}). Then for  $1< p< \infty$,
\end{theorem}
\begin{equation}\label{ineq:MzInfty}
\lim\sup_{n\rightarrow\infty} \max_{P_n(z)\in\Pi_n} \frac{\int_{\gamma_0} |P_n(z)|^p |dz| }{\sum^n_{k = 0} |P_n(z^*_{n,k})|^p \text{dist}(z^*_{n,k},\Gamma_n)} = \infty.
\end{equation}

As a consequence of Theorem 2, the ``modified" Marcinkiewicz-Zygmund inequalities defined in (\ref{mod MZ ineq}) fail for some polynomials $P_n(z)\in \Pi_n$ with $n=n_j\rightarrow \infty$. To establish the third theorem of this paper, we will describe a method for adjusting the Fej\'er points ${\bf z}^{*}_n =\{z^*_{n,k}\}_{k=0}^n \subset \gamma_0$ in the next section, and let ${\bf{z}^{**}_n} = \{ z^{**}_{n,k}\}_{k=0}^n$ denote the adjustment of ${\bf{z}^{*}_n}$.

\begin{theorem}  \label{theorem:AdJustFejer}
For the $L$-shape arc $\gamma_0\subset \CC$ defined in $(\ref{definegamma})$, there exists a family of adjusted Fej\'er points  ${\bf{z}^{**}_n} = \{ z^{**}_{n,k}\}_{k=0}^n\subset\gamma_0$ and two  constants $c_1>0$ and $c_2 > 0$, such that  
\begin{equation} \label{ineq:Lebesg2}
c_1 \log^2(n) \le  L_{{\bf {z^{**}_n}}}  \le c_2 \log^2(n),
\end{equation}
for all $n\in \ZZ_{+}$.
\end{theorem}

Throughout the presentation in the next three sections, the positive constants $c_1$ and $c_2$ are ``varying" constants in that they may be different at different occasions, even within the formulation of the same formula.\\ 

This paper is organized as follows. In the next section, we will discuss the distribution of Fej\'er points on the $L$-shape open arc $\gamma_0\subset \CC$ defined in (\ref{definegamma}) and describe an adjustment of these points. Preliminary lemmas, particularly on the properties of $\omega_n(z)$ defined in (\ref{Omega0}), will be established in Section 3. Theorems 1--3 will be proved in Section 4. The final section of the paper will be devoted to the discussion of the notion of Lebesgue constants for the general ``polynomial preservation" linear operators $A_n$ on Banach spaces, the limitations of Fej\'er points  on open arcs in the complex plane, approximation by polynomials with restricted zeros and the related results on asymtoptically neutral distribution of electrons, based on and inspired by the work of J. Korevaar in the decade of the 1960's, as well as the problem of placement of electrons on the surfaces on bounded and connected domains in $\RR^3$, governed by the ``total energy" due to the locations of the electrons; and this leads to the discussion of the inverse problem of super-resolution recovery of point-masses in $\RR^s$ for any $s\ge 1$.\\ 

In the Epilogue of this paper, we point out the ``potential" research opportunities of the topic of point-mass representation, understanding, manipulation, and recovery, particularly for the advancement of {\bf S}ciences, {\bf T}echnology and {\bf E}ngineering, facilitated by the advancement of {\bf M}athematics: leading to a productive and mutually beneficial cross-fertilization platform, commonly called ``{\bf STEM}".
\vskip .1in

\bhag{Distribution of Fej\'er points and their adjustment}

Since the extension to the boundary of the conform map $\psi_0$ is not a one-to-one mapping of the unit circle of $|w| = 1$ onto the open arc $\gamma_0$ (with the exception of the 2 end-points $\psi_0(e^{\pm i{2\pi/3}})$ of $\gamma_0$), in that two points on $|w|=1$ could be mapped to a single interior point of $\gamma_0$, but the relation of these two points on $|w|=1$ can be adjusted. To do so, let $J(t): [{{2\pi}\over{3}}, \pi] \rightarrow  [0, {{2\pi}\over{3}}]$ be defined to satisfy the condition:
\begin{equation}
\psi_0(e^{i t}) = \psi_0(e^{iJ(t)}),
\end{equation}
and for $t\in [-\pi, -{{2\pi}\over{3}}]$, we set $J(t) = -J(-t)$.
Then it follows from ($\ref{ineq:Conform0}$) that 
\begin{equation} \label{ineq:sin  J(t)}
\sin(J(t))\sin^2(  J(t)/2) = \sin(t)\sin^2(t/2), 
\end{equation}
so that the formula of $J(t)$ for $t\in [{{2\pi}\over{3}},\pi]$ is given by
$$J(t) = 2\arcsin\big(\sqrt{( (1-\sin^2(t/2)) ( 1 + (1 + 2\sin^2(t/2)) / r(t) ) + r(t))/3}\,\,\big),$$
with
$$
r(t) =  \sqrt[3]{(1-\sin^2(t/2))(2 + 5\sin^2(t/2) + 20\sin^4(t/2) + \sin^2(t/2)\sqrt{27(3+8\sin^2(t/2)+16\sin^4(t/2))})/2}.
$$
Although the formula is seemingly quite complicated, it is precise with the only ``singularity" at $t = \pi$. From ($\ref{ineq:sin  J(t)}$), we  also have
\begin{equation} \label{ineq:J(t)ClosePi}
  J(t) = \sqrt[3]{4}(\pi - t)^{1/3} + (\pi - t)/3 + O((\pi - t)^{5/3}),
\end{equation}
and by taking the derivative on both sides of ($\ref{ineq:sin  J(t)}$), we also obtain 
$$
J'(t) = {{\cos(t) - \cos(2t)}\over{\cos(  J(t)) - \cos(2  J(t))}},
$$
so that it follows from ($\ref{ineq:J(t)ClosePi}$) that
\begin{equation} \label{ineq:J'(t)ClosePi}
J'(t) = -\sqrt[3]{4}/3(\pi - t)^{-2/3}  - 1/3 + o(1),
\end{equation}
for $t\rightarrow{\pi}$. Note that both ($\ref{ineq:J(t)ClosePi}$) and ($\ref{ineq:J'(t)ClosePi}$) are mainly derived for the values of $t$ that are close to $\pi$. In this regard, note that from ($\ref{ineq:J'(t)ClosePi}$), we have $J'(t) < -1$ when $t$ is close to $\pi$. The solution of the equation $J'(t) = -1$, under the constraint of ($\ref{ineq:sin  J(t)}$), is $t_0 = \pm {{2\pi}\over{3}}$. On the other hand, since $\lim_{t\rightarrow{\pi}} J'(t) = -\infty$, it follows that

\begin{equation} \label{ineq:J'(t)Ge1}
J'(t) < -1, \hskip 15pt {{2\pi}\over{3}}< |t| < \pi.
\end{equation}
Observe that since $\psi_0(e^{it}) = \psi_0(e^{i  J(t)})$, there might be two distinct points $\theta_{n,j}$ and $J(\theta_{n,k})$ that are very close to each other. For this situation, we need to establish the following result on the existence of some upper bound.

\begin{lemma}\label{lemma:closeToJ_theta}
There exists a constant $c_2> 0$, such that for every positive integers $j$, 
\begin{equation}\label{ineq:closeToJ_theta}
 \min_{n\ge j}(n+1)^{\frac{4}{3}}|J(\theta_{n,{\lfloor{n/2}\rfloor}} )- \theta_{n,j}| \le c_2.
\end{equation}
\end{lemma}  

\begin{Proof}
Since $\theta_{n,{\lfloor{n/2}\rfloor}} = \pi - \frac{\pi}{n+1}$, it follows from ($\ref{ineq:J(t)ClosePi}$) that
\begin{equation}\label{ineq:thetaHalfN}
r_n := (n+1) J(\theta_{n,{\lfloor{n/2}\rfloor}} ) = \sqrt[3]{4\pi}(n+1)^{2/3} + \frac{\pi}{3} + O((\frac{\pi}{n+1})^{2/3}),
\end{equation} 
so that 
\begin{equation}\label{ineq:thetaHalfDiff}
 r_{n+2} - r_{n} = (n+3)J(\theta_{n+2,{\lfloor{n/2+1}\rfloor}} ) - (n+1)J(\theta_{n,{\lfloor{n/2}\rfloor}} )  = \frac{4}{3}\sqrt[3]{4}(\frac{\pi}{n+1})^{1/3} +  O((\frac{\pi}{n+1})^{2/3}). \nonumber
\end{equation}
Hence, there exists some positive integer $n_0$, such that for $n \ge n_0$,
\begin{equation}\label{ineq:thetaHalfDiff1}
 \frac{c_1}{\sqrt[3]{n+1}} \le r_{n+2} - r_{n} \le \frac{c_2}{\sqrt[3]{n+1}},
\end{equation}
and $\lim_{n\rightarrow\infty} r_{2n} = \infty$. 

Next, for any integer $j > \frac{r_{n_0}}{2\pi}$, let $n(j)$ be an even integer, for which $r_{n(j)} \le 2\pi j \le r_{n(j) + 2}$. From $(\ref{ineq:thetaHalfDiff1})$, we have
$$
0 \le 2\pi j - r_{n(j)} < \frac{c_2}{\sqrt[3]{n(j)+1}}.
$$
As a consequence, 
$$
({n(j)+1})^{\frac{4}{3}}|\theta_{n(j), j} - \theta_{n(j), \frac{n(j)}{2}}| \le c_2, \hskip 10pt j > \frac{r_{n_0}}{2\pi},
$$
so that $(\ref{ineq:closeToJ_theta})$ holds for $j > \frac{r_{n_0}}{2\pi}$. Similarly, we can show that it also holds for all $j > 0$ with sufficiently large $c_2$.
\end{Proof}

One the other hand, it was already shown in $\cite{CZ2021}$ that for the validity of the modified Marcinkiewicz-Zygmund inequalities defined in (\ref{mod MZ ineq}), the interpolation nodes must satisfy the separation condition: 

\begin{equation}\label{sep cond1}
\min_{|\theta_{n,j}|<{{2\pi}\over{3}}, |\theta_{n,k}|\ge {{2\pi}\over{3}} }  |\theta_{n,j} - J(\theta_{n,k})| \ge  {{c_1}\over{n+1}}. 
\end{equation}
By applying Lemma $\ref{lemma:closeToJ_theta}$,  we will prove in Section 4 that the modified Marcinkiewicz-Zygmund inequalities are invalid for the family ${\bf z}^{*}_n$ of Fej\'er points with $n=n(j)$. 
\vskip .1in
In the following, we describe a method for adjusting the Fej\'er points on $\gamma_0$ to meet the separation requirement (\ref{sep cond1}). To do so, let us first describe the method for adjusting those Fej\'er points that correspond to $0 \le k \le \lfloor{n/2}\rfloor$; that is, $0 < \theta_{n,k} < \pi$ by applying the definition of $\theta_{n,k}$ in $(\ref{ineq:theta0})$.
\vskip .1in

{\bf Introducing Adjusted Fej\'er points on the open arc $\gamma_0$}
\vskip .1in
(1) For any $0\le \theta_{n,j} \le \frac{2\pi}{3}$, we introduce
\begin{equation}\label{ineq:thetaTildeJ}
\tilde\theta^*_{n,j} = \left\{
\begin{array}{ll}
      \theta_{n,j}, & \min_{\{k:\theta_{n,k}>2\pi/3\}}|\theta_{n,j} - J(\theta_{n,k})| \ge \frac{2\pi}{3(n+1)}; \\
        J(\theta_{n,k}) + \frac{2\pi}{3(n+1)}, & \text{if there exists $k$ such that } \theta_{n,j} - \frac{2\pi}{3(n+1)} < J(\theta_{n,k}) \le \theta_{n,j};\\
            J(\theta_{n,k}) - \frac{2\pi}{3(n+1)}, & \text{if there exists $k$ such that } \theta_{n,j}  < J(\theta_{n,k}) < \theta_{n,j} + \frac{2\pi}{3(n+1)}. \\
\end{array} 
\right. 
\end{equation}

(2) For any $\frac{2\pi}{3} < \theta_{n,k} < \pi$, we introduce
\begin{equation}\label{ineq:thetaTildeK}
\tilde\theta^*_{n,k} = \left\{
\begin{array}{ll}
      J(\theta_{n,k}), & \min_{\{j:\theta_{n,j}\le 2\pi/3\}}|\theta_{n,j} - J(\theta_{n,k})| \ge \frac{2\pi}{3(n+1)};  \\
       \theta_{n,j} - \frac{2\pi}{3(n+1)},  & \text{if there exists $j$ such that } \theta_{n,j} - \frac{2\pi}{3(n+1)} < J(\theta_{n,k}) \le \theta_{n,j};\\
             \theta_{n,j} + \frac{2\pi}{3(n+1)}, & \text{if there exists $j$ such that } \theta_{n,j}  < J(\theta_{n,k}) < \theta_{n,j} + \frac{2\pi}{3(n+1)}. \\
\end{array} 
\right. 
\end{equation}

(3) For $\pi\le \theta_{n,k}< 2\pi$; that is, for $\lfloor{n/2}\rfloor < k \le n$, we define $\tilde\theta^*_{n,k}$ by anti-symmetry with respect to $\pi$, namely:
\begin{equation}\label{ineq:tildaJKNeg}
\tilde\theta^*_{n,k} = -\tilde\theta^*_{n,2\lfloor{n/2}\rfloor + 1 - k}, \hskip 10pt \lfloor{n/2}\rfloor < k \le n.\\
\end{equation}

(4) The definition of $\{\tilde\theta^*_{n,k}\}_{k=0}^n$ from the totality of (\ref{ineq:thetaTildeJ}), (\ref{ineq:thetaTildeK}) and (\ref{ineq:tildaJKNeg}) are applied to arrive at the following family of adjusted Fej\,er points ${\bf z}^{**} :=\{\z^{**}_{n,k}\}_{k=0}^n$ that lie on the open arc $\gamma_0$, defined by:
\begin{equation}\label{adjFejerPoints}
z^{**}_{n,k} = \psi_0(\tilde\theta^*_{n,k}), \,\,\,\,\,k = 0,\cdots,n.
\end{equation}
{\bf End of definition of adjusted Fej\'er points.}\,\,This is the family of adjusted Fej\'er points in the statement of Theorem $\ref{theorem:AdJustFejer}$ in the first section.  
\vskip .1in
Under the condition of definition $(\ref{ineq:thetaTildeJ})$ and  (\ref{ineq:thetaTildeK}), we have
 \begin{equation} \label{adjDiffJK}
 \left\{
\begin{array}{ll}
      |\theta_{n,j} -\tilde\theta^*_{n,j}| \le \frac{2\pi}{3(n+1)}, & |\theta_{n,j}| \le \frac{2\pi}{3}; \\

    |J(\theta_{n,k}) -\tilde\theta^*_{n,k}| \le \frac{2\pi}{3(n+1)}, & |\theta_{n,k}| > \frac{2\pi}{3};\\
\end{array} 
\right. 
\end{equation}
and under the conditions $|\theta_{n,j}| \le \frac{2\pi}{3}$ and $|\theta_{n,k}| > \frac{2\pi}{3}$, if, in addition, $|\theta_{n,j} - J(\theta_{n,k})| < \frac{2\pi}{3(n+1)}$, then we will call $(j, k)_n$ an ``adjustment pair",\,\,and remark that only adjustment pairs require adjustment. 

The objective of the following lemma is for establishing Theorem $\ref{theorem:AdJustFejer}$ in Section 4 and the proof of its preliminary lemmas in Section 3.

\begin{lemma}\label{lemma:separate_theta}
For the $\{\tilde\theta^*_{n,k}\}$ defined in $(\ref{ineq:thetaTildeJ})$ and $(\ref{ineq:thetaTildeK})$, 
\begin{equation} \label{ineq:separate0}
\min_{j\neq k} |\tilde\theta^*_{n,j} - \tilde\theta^*_{n,k}| \ge {{2\pi}\over{3(n+1)}}.
\end{equation}
\end{lemma}  
\begin{Proof}
Since there are several cases in the description of the adjustment $\{\tilde\theta^*_{n,k}\}$, we will verify $(\ref{ineq:separate0})$ in $6$ cases, as follows:
\vskip.1 in
\noindent case 1. Both $\tilde\theta^*_{n,j}$ and $\tilde\theta^*_{n,k}$ are defined by $(\ref{ineq:thetaTildeJ})$;

\noindent case 2. Both $\tilde\theta^*_{n,j}$ and $\tilde\theta^*_{n,k}$ are defined by $(\ref{ineq:thetaTildeK})$; 

\noindent case 3.  $\tilde\theta^*_{n,j}$ and $\tilde\theta^*_{n,k}$ are defined by  $(\ref{ineq:thetaTildeJ})$ and $(\ref{ineq:thetaTildeK})$ respectively. $(j,k)_n$ is an adjustment pair;

\noindent case 4.  $\tilde\theta^*_{n,j}$ and $\tilde\theta^*_{n,k}$ are defined by  $(\ref{ineq:thetaTildeJ})$ and $(\ref{ineq:thetaTildeK})$ respectively. neither $j$ nor $k$ is in any adjustment pair;

\noindent case 5.  $\tilde\theta^*_{n,j}$ and $\tilde\theta^*_{n,k}$ are defined by  $(\ref{ineq:thetaTildeJ})$ and $(\ref{ineq:thetaTildeK})$ respectively.  $j$ is in an adjustment pair $(j,k_1)_n$ but $k \neq k_1$;

\noindent case 6.  $\tilde\theta^*_{n,j}$ and $\tilde\theta^*_{n,k}$ are defined by  $(\ref{ineq:thetaTildeJ})$ and $(\ref{ineq:thetaTildeK})$ respectively.  $k$ is in an adjustment pair $(j_1,k)_n$ but $j \neq j_1$.\\

\noindent For case 1, it follows from $(\ref{adjDiffJK})$ that
$$
 |\tilde\theta^*_{n,j} - \tilde\theta^*_{n,k}| =  |\theta_{n,j} - \theta_{n,k} + \tilde\theta^*_{n,j} - \theta_{n,j} + \theta_{n,k} - \tilde\theta^*_{n,k}| \ge \frac{2\pi}{n+1} - \frac{2\pi}{3(n+1)} - \frac{2\pi}{3(n+1)}.
$$

\noindent For case 2, combining with $(\ref{adjDiffJK})$ and $|J'(t)| > 1$, we have
$$
|\tilde\theta^*_{n,j} - \tilde\theta^*_{n,k}| \ge |J(\theta_{n,j}) - J(\theta_{n,k})| - |\tilde\theta^*_{n,j} - J(\theta_{n,j})| - |J(\theta_{n,k}) - \tilde\theta^*_{n,k}|\ge \frac{2\pi}{n+1} - \frac{2\pi}{3(n+1)} - \frac{2\pi}{3(n+1)}. 
$$

\noindent For case 3, if $(j,k)_n$ is an adjustment pair, it follows from $(\ref{ineq:thetaTildeJ})$ and $(\ref{ineq:thetaTildeK})$ that 
$$
|\tilde\theta^*_{n,j} - \tilde\theta^*_{n,k}| = \frac{4\pi}{3(n+1)} - |\theta_{n,j} - J(\theta_{n,k})| \ge  \frac{2\pi}{3(n+1)}. 
$$

\noindent For case 4, as neither $j$ nor $k$ are in any adjustment pair, we have
$$
|\tilde\theta^*_{n,j} - \tilde\theta^*_{n,k}| = |\theta_{n,j} - J(\theta_{n,k})| \ge \frac{2\pi}{3(n+1)}.
$$

\noindent For case 5, by combining with $(\ref{ineq:thetaTildeJ})$  and $(\ref{adjDiffJK})$, we obtain
$$
|\tilde\theta^*_{n,j} - \tilde\theta^*_{n,k}| \ge | J(\theta_{n,k_1}) - J(\theta_{n,k})| - | J(\theta_{n,k_1}) - \tilde\theta^*_{n,j}| - |J(\theta_{n,k}) - \tilde\theta^*_{n,k}|  \ge \frac{2\pi}{3(n+1)}.
$$

\noindent For case 6, by combining with $(\ref{ineq:thetaTildeK})$  and $(\ref{adjDiffJK})$, we obtain
$$
|\tilde\theta^*_{n,j} - \tilde\theta^*_{n,k}| \ge | \theta_{n,j} - \theta_{n,j_1}| - | \tilde\theta^*_{n,j} - \theta_{n,j}| - |\theta_{n,j_1} - \tilde\theta^*_{n,k}|  \ge \frac{2\pi}{3(n+1)}.
$$

\noindent This completes the proof of the lemma. 
\end{Proof}

\bhag{Preliminary results}
The derivation of our preliminary lemmas in this section depends heavily on the level curves: 
\begin{equation} \label{level curves}
\begin{cases}
 {\displaystyle \rho_m:= 1 + \frac{1}{m+1}}; \\
 &\\ 
{\displaystyle \Gamma_m := \psi_0 (|w| = \rho_m)},\\   
\end{cases}
\end{equation} 
for integers $1<m\in \ZZ$, where the level curves $\Gamma_m $ enclose the open arc $\gamma_0$ and converge to $\gamma_0$ as $m\rightarrow \infty$.

\begin{lemma}\label{lemma:omega_bound2}
Let $\{\zeta^*_{n,k} = \psi_0\big(\rho_n e^{i\theta_{n,k}}\big)\}_{k = 0}^n$ and $\omega^*_n(z) = \prod_{k = 0}^n(z - \zeta^*_{n,k})$. Then
\begin{equation} \label{ineq:omega_bound2}
e^{-3}(e - 1)^2 \le |\omega^*_n(z)| \le \frac{e^{5}(1 + 2 e)}{e - 1}, \hskip 10pt z\in \gamma_0.
\end{equation} 
\end{lemma}        

The proof of this lemma is similar to that of Lemma 2 in $\cite{Z1995}$, which depends on the property of  bounded rotation in Lemma 3 of the paper \cite{KOV1966} on approximation by Faber polynomials.
\vskip 15pt
\begin{Proof}
As in $\cite{CUR1935}$, for $|u| \ge 1$ and $|w| \ge 1$, we define
$$
\Psi(w, u) = \left\{
\begin{array}{ll}
      \log\big({{\psi_0(w) - \psi_0(u)}\over{w - u}});    &u \neq w, \\
      \\
       \log(\psi_0'(w)),      &u = w,\\
\end{array} 
\right. 
$$
so that
$$
\Psi(w, u) = \sum_{m = 1}^{\infty} {{a_m(w)}\over{u^m}}.
$$
For $|w| \ge 1$,
$$
a_m(w) = {1\over{2\pi i}}\int_{|u| = |w| + \epsilon} u^{m-1}\Psi(w, u) du = {1\over{2\pi m i}}\int_{|u| = |w| + \epsilon} u^{m}\big({1\over{u - w}} -{{\psi_0'(u)}\over{\psi_0(u) - \psi_0(w)}}\big) du = {{w^m - F_m(\psi_0(w))}\over{m}},
$$
where $F_m(z)$ is the $m$-th degree Faber polynomial respect to $\gamma_0$ (resulting from those corresponding $m$-th degree Faber polynomials defined in the open region $\CC\backslash \bar D_N$, by taking $N\rightarrow \infty$, where $D_N$ is the bounded domain with $\partial D_N = \Gamma_N$). From the integral formula of Faber polynomials in $\cite{KOV1966}$,
$$
F_m(\psi_0(e^{it})) = {1\over\pi}\int_0^{2\pi}e^{ims}d_s(\arg(\psi_0(e^{is})-\psi_0(e^{it}))),
$$
we have
\begin{equation} \label{a_mW}
|a_m(w)|\le {{4}\over{m}}, \hskip 10pt |w| = 1.    
\end{equation}
For $z = \psi_0(w) \in\gamma_0$, we have
$$
\log\big({{\omega_n^*(z)}\over{w^{n+1} - (-1)^{n}\rho_n^{n+1}}}\big) = \sum_{k=0}^n\Psi(w,\rho_n e^{\theta_{n,k}}) = \sum_{m = 1}^{\infty}\rho_n^{-m}a_m(w) \sum_{k=0}^n e^{-i m\theta_{n,k}}.
$$
On the other hand, since
$$
\sum_{k=0}^n e^{-i m\theta_{n,k}} = \left\{
\begin{array}{ll}
      (-1)^{nj}(n+1),    &m = (n+1)j \\
      \\
       0,          &\,\, \text{otherwise} \,,\\
\end{array} 
\right. 
$$
it follows from $(\ref{a_mW})$ that
$$
\big|\log\big({{\omega_n^*(z)}\over{w^{n+1} - (-1)^{n}\rho_n^{n+1}}}\big) \big| = \big|(n+1)\sum_{j = 1}^{\infty} {{(-1)^{nj}a_{(n+1)j}(w)}\over {\rho_n ^{(n+1)j}}}\big| \le 4\sum_{j = 1} ^{\infty} {1\over{j e^j}} = -4\log(1-e^{-1}),
$$
so that 
$$
|\omega_n^*(z)| \le {e^{5}\over{e - 1}} |{w^{n+1} - (-1)^{n}\rho_n^{n+1}}| \le e^{5}{(1 + 2 e)\over{e - 1}}, \hskip 10pt z\in \gamma_0,
$$
and
$$
|\omega_n^*(z)| \ge e^{-3}(e - 1) |{w^{n+1} - (-1)^{n}\rho_n^{n+1}}| \ge e^{-3}(e - 1)^2 , \hskip 10pt z\in \gamma_0.
$$
This completes the proof of the lemma.
\end{Proof}

While the points  $\zeta^*_{n,k}$, $k=0,\cdots,n$, do not lie on the open arc $\gamma_0$, our goal is to replace the Fej\'er points $\{z^*_{n,k}\}_{k=0}^n\subset \gamma_0$ by $\{\zeta^*_{n,k}\}_{k=0}^n\subset\Gamma_n$. In the following,  we will establish an approximation result for such replacement.

For any $t\in [0, \frac{2\pi}{3})$ and $t\in[\frac{2\pi}{3},\pi)$, let $k_1(t)$ and $k_2(t)$ be defined, respectively, as follows:
\begin{equation} \label{k1k2}
\begin{cases}
{\displaystyle k_1(t) := k_{1,n}(t) = \arg\min_{\{k:0\le\theta_{n,k}\le \frac{2\pi}{3}\}}\big|\theta_{n, k} - t\big|};\\
&\\ 
    {\displaystyle k_2(t) := k_{2,n}(t) =\arg\min_{\{k:\frac{2\pi}{3} < \theta_{n,k} < \pi\}}\big|\theta_{n, k} - J^{-1}(t)\big|}.
\end{cases}
\end{equation} 
 Also, for $t\in [-{{2\pi}\over{3}},0)$ and $t\in [-\pi, -\frac{2\pi}{3})$, we set $k_1(t) = \lfloor{n/2}\rfloor + 1 - k_1(-t)$ and $k_2(t) =  \lfloor{n/2}\rfloor + 1 - k_2(-t)$, respectively, by applying anti-symmetric extension with respect to $0$. For each of the above four cases, if the minimum is attained at two points $\theta_{n,k}$, we choose the smaller index $k$.

In this paper, we say that two quantities $A(\cdots)$ and $B(\cdots)$ are equivalent, if there exist two constants $c_2 \ge c_1 > 0$, such that
$$
c_1 A(\cdots) \le B(\cdots) \le c_2 A(\cdots),
$$
and introduce the notion
$$
 A(\cdots) \tilde{=} B(\cdots).
$$

We need the separation of $u$ and $w$ under the conformal mapping. From \cite{D1977} (see  page 387),  we obtain, for $1\le |u|, |w|\le 2$, the following equivalence relationship:
\begin{equation} \label{ineq:dist_00}
{{|\psi_0(u) - \psi_0(w)|}} \tilde{=} |u-w|(|u - 1|+|u-w|)^{{1\over{2}}},|\arg(u)|\in[0, \pi/3],  |\arg(w)| \in [0, {{2\pi}\over{3}}], 
\end{equation}

\begin{equation} \label{ineq:dist_01}
{{|\psi_0(u) - \psi_0(w)|}} \tilde{=} |u-w|(|u - e^{i {2\over 3} \pi}|+|u-w|),\arg(u)\in[\pi/3, {{2\pi}\over{3}}],  \arg(w) \in [0, {{2\pi}\over{3}}],
\end{equation}

\begin{equation} \label{ineq:dist_10}
{{|\psi_0(u) - \psi_0(w)|}} \tilde{=} |u-w|(|u - e^{i {2\over 3} \pi}|+|u-w|), \arg(u)\in[{{2\pi}\over{3}}, 5\pi/6],  |\arg(w)| \in [{{2\pi}\over{3}}, \pi], 
\end{equation}

\begin{equation} \label{ineq:dist_11}
{{|\psi_0(u) - \psi_0(w)|}} \tilde{=} |u-w|(|u + 1|+|u-w|)^{-{1\over{2}}},|\arg(u)|\in[5\pi/6, \pi],  |\arg(w)| \in [{{2\pi}\over{3}}, \pi], 
\end{equation}
and
\begin{equation} \label{ineq:dist_12}
{{|\psi_0(u) - \psi_0(w)|}} \tilde{=} |u-w|(|u - e^{-i {2\over 3} \pi}|+|u-w|),\arg(u)\in[-5\pi/6, -{{2\pi}\over{3}}],  |\arg(w)| \in [{{2\pi}\over{3}}, \pi]. 
\end{equation}

Since estimation of the distances between $u$ and $w$ under the conformal mapping is very tedious, we partition the interval $[0,\pi]$ into 4 sub-intervals: $I_1 = [0, \pi/3), I_2 = [\pi/3, {{2\pi}\over{3}}]$, $I_3 = ({{2\pi}\over{3}}, 5\pi/6]$ and $I_4 = (5\pi/6, \pi]$, and will show that each sub-interval $I_j$, for $j = 1,2,3,4,$ contains exactly one singular end-point $s_j$, in that $|\psi_0'(e^{i s_j})| = 0$ or $\infty$. The same argument also applies to the interval $[-\pi,0]$.

If $\arg(w)$ and $\arg(u)$ are in the same sub-interval, say $I_j$, let ${v}\in \{u,w\}$, such that $|{v} - e^{i s_j}| = \max\{|u - e^{i s_j}|, |w - e^{i s_j}|\}$. Since $|{v} - e^{i s_j}|\tilde{=} (|u - e^{i s_j}| +|w - u|)$, it follows from ($\ref {ineq:dist_00})$ and $(\ref{ineq:dist_12}$) that if $\arg(u)$ and $\arg(w)$ are in the same sub-interval $I_j$, then
\begin{equation} \label{ineq:dist21}
{{|\psi_0(u) - \psi_0(w)|}} \tilde{=}  |\psi_0'({v})| |u-w|. 
\end{equation}
On the other hand, if $\arg(u)\in I_1$ and $\arg(w)\in I_2$, or if $\arg(u)\in I_3$ and $\arg(w)\in I_4$, then
\begin{equation} \label{ineq:dist22}
{{|\psi_0(u) - \psi_0(w)|}} \tilde{=} |u-w|. 
\end{equation}

\begin{lemma} \label{lemma:Omega_ToOmega1}
Let $z = \psi_0(e^{i t})\in\gamma_0$ with $|t|\in [0, 2/3\pi]$. Then for $0\le |t|\le{{2\pi}\over3}$,
\begin{equation}\label{ineq:Key0}
|\omega_n(z)| \tilde{=} \big |{{ (z - z^*_{n,k_1(t)}) ( z - z^*_{n,k_2(t)}) } \over {(z - \zeta^*_{n,k_1(t)}) ( z - \zeta^*_{n,k_2(t)})}}\big|;
\end{equation}
and for $|\theta_{n,j}| \in [0,2/3\pi]$,
\begin{equation}\label{ineq:Key'01}
|\omega'_n(z^*_{n,j})| \tilde{=} \big |{{ z^*_{n,j} - z^*_{n,k_2(\theta_{n,j})} } \over {(z^*_{n,j} - \zeta^*_{n,j}) ( z^*_{n,j} - \zeta^*_{n,k_2(\theta_{n,j})})}}\big|.
\end{equation}
\end{lemma}        

\begin{Proof} From ($\ref{ineq:omega_bound2}$) in Lemma $\ref{lemma:omega_bound2}$, it is sufficient to prove
\begin{equation} \label{ineq:Omega_ToOmega1}
c_1\le \big|{{\omega^*_n(z)/( (z - \zeta^*_{n,k_1(t)}) ( z - \zeta^*_{n,k_2(t)}) )}\over{{\omega_n(z)/( (z - z^*_{n,k_1(t)}) ( z - z^*_{n,k_2(t)}) )}   }}\big|\le c_2,
\end{equation}
for $t \in [0, {{2\pi}\over{3}}]$, since the proof for $t \in [-{{2\pi}\over{3}}, 0)$ is only a consequence of the symmetry of $\gamma_0$. For $t \in [0, {{2\pi}\over{3}}]$, it is obvious that
\begin{equation}\label{eq:sep_sum}
\big |{{z - \zeta^*_{n,k}}\over{z - z^*_{n,k}}}\big|^2 = 1 + \big |{{z^*_{n,k} - \zeta^*_{n,k}}\over{z - z^*_{n,k}}}\big|^2 + 2\operatorname{Re}\big({{z^*_{n,k} - \zeta^*_{n,k}}\over{z - z^*_{n,k}}}\big).
\end{equation}
Under the condition of $k \neq k_1(t)$ and $k\neq k_2(t)$, it is also clear that $|e^{it} - e^{i \theta_{n,k}}| \tilde{=}  |e^{it} - \rho_n e^{i \theta_{n,k}}|$. From the distance estimation $(\ref{ineq:dist_00})$,  we have, for $\theta_{n,k} \in I_1$,
$$
|z - z^*_{n,k}| \tilde{=} |e^{it} - e^{i \theta_{n,k}}|\sqrt{|1 -  e^{i \theta_{n,k}}| + |e^{it} - e^{i \theta_{n,k}}|}\tilde{=} |e^{it} - \rho_n e^{i \theta_{n,k}}|\sqrt{|1 - \rho_n e^{i \theta_{n,k}}| + |e^{it} - \rho_n e^{i \theta_{n,k}}|}\tilde{=} |{z - \zeta^*_{n,k}}|.
$$
Similarly, as long as $k \neq k_1(t)$ and $k\neq k_2(t)$, we have
$$
|z - z^*_{n,k}| \tilde{=} |z - \zeta^*_{n,k}|,
$$
which is equivalent to
\begin{equation} \label{log2Diff1}
\big|\log \big |{{z - \zeta^*_{n,k}}\over{z - z^*_{n,k}}}\big|^2 \big| \tilde{=} \big|\big |{{z - \zeta^*_{n,k}}\over{z - z^*_{n,k}}}\big|^2 - 1\big |, \hskip 10pt k\neq k_1(t), k_2(t).
\end{equation}
Therefore, it follows from $(\ref{eq:sep_sum})$ that
\begin{equation} \label{sumLog2Diff1}
\big|\log\big|{{\omega^*_n(z)/( (z - \zeta^*_{n,k_1(t)}) ( z - \zeta^*_{n,k_2(t)}) )}\over{{\omega_n(z)/( (z - z^*_{n,k_1(t)}) ( z - z^*_{n,k_2(t)}) )}   }}\big|^2\big | \le c_2\sum_{k\neq k_1(t), k_2(t)}\big(\big |{{z^*_{n,k} - \zeta^*_{n,k}}\over{z - z^*_{n,k}}}\big|^2 + 2\big|\operatorname{Re}\big({{z^*_{n,k} - \zeta^*_{n,k}}\over{z - z^*_{n,k}}}\big)\big|\big).
\end{equation}
The proof of (\ref{ineq:Omega_ToOmega1}) will be complete by showing that the left-hand side of above inequality is bounded. In this regard, let us consider the sum of the ``square terms" in the left-hand side of $(\ref{sumLog2Diff1})$, and observe that from (\ref{ineq:dist21}),
\begin{equation}
\sum_{k\neq k_1(t), k_2(t)}\big |{{z^*_{n,k} - \zeta^*_{n,k}}\over{z - z^*_{n,k}}}\big|^2 
\le   c_2\sum_{k\neq k_1(t), k_2(t)} {{n^{-2}|\psi_0'(\rho_n e^{i\theta_{n,k}})|^2}\over{|z - z^*_{n,k}|^2}}.
\end{equation}
For $k\neq k_1(t)$, it is obvious that $|e^{i t}-e^{i\theta_{n,k}}|\ge\pi/(n+1)$ for $\theta_{n,k} \in [0, {{2\pi}\over{3}}]$, and under the condition of $\theta_{n,k}\in I_1$, we have
\begin{equation} \label{ineq:sq10}
|\psi_0(e^{i \theta_{n,k}}) - \psi_0(e^{it})|\ge c_1|e^{i \theta_{n,k}} - e^{it}|(|e^{i \theta_{n,k}} - 1| + \pi/(n+1))^{1/2} \tilde{=} |e^{i \theta_{n,k}} - e^{it}| |\psi_0'(\rho_n e^{i\theta_{n,k}})|.
\end{equation}
Similarly, for $\theta_{n,k}\in I_2$, we have
\begin{equation} \label{ineq:sq11}
|\psi_0(e^{i \theta_{n,k}}) - \psi_0(e^{it})|\ge c_1|e^{i \theta_{n,k}} - e^{it}|(|e^{i \theta_{n,k}} - e^{i{{2\pi}\over{3}}}| + \pi/(n+1)) \tilde{=} |e^{i \theta_{n,k}} - e^{it}| |\psi_0'(\rho_n e^{i\theta_{n,k}})|.
\end{equation}
In summary, for $\theta_{n,k} \in [0, {{2\pi}\over{3}}]$,
\begin{equation} \label{ineq:sq12}
|z - z^*_{n,k}| = |\psi_0(e^{i \theta_{n,k}}) - \psi_0(e^{it})|\ge c_1 |e^{i \theta_{n,k}} - e^{it}| |\psi_0'(\rho_n e^{i\theta_{n,k}})|.
\end{equation}
In addition, it follows from $(\ref{ineq:sq10})$ that
\begin{equation} \label{ineq:sum_sq00}
\sum_{k\neq k_1(t), \theta_{n,k} \in [0, {{2\pi}\over{3}}]}\big |{{z^*_{n,k} - \zeta^*_{n,k}}\over{z - z^*_{n,k}}}\big|^2 
\le  c_2\sum_{k\neq k_1(t), \theta_{n,k} \in [0, {{2\pi}\over{3}}]} {{n^{-2}}\over{|e^{i t} - e^{i \theta_{n,k}}|^2}} 
\le c_2.
\end{equation}

Next, for $\theta_{n,k}\in I_3$,
\begin{equation} \label{ineq:sq20}
|z - z^*_{n,k}| = |\psi_0(e^{i \theta_{n,k}}) - \psi_0(e^{iJ^{-1}(t)})|\ge c_1 |e^{i \theta_{n,k}} - e^{iJ^{-1}(t)}| |\psi_0'(\rho_n e^{i\theta_{n,k}})|.
\end{equation}
Similarly to $(\ref{ineq:sum_sq00})$,
\begin{equation} \label{ineq:sum_sq10}
\sum_{k\neq k_2(t),\theta_{n,k}\in I_3}\big |{{z^*_{n,k} - \zeta^*_{n,k}}\over{z - z^*_{n,k}}}\big|^2 
\le  c_2\sum_{k\neq k_2(t)} {{n^{-2}}\over{|e^{iJ^{-1}(t)} - e^{i \theta_{n,k}}|^2}} 
\le c_2.
\end{equation}

For $\theta_{n,k}\in I_4$, it follows from ($\ref{ineq:dist_11}$) that
$$
|z - z^*_{n,k}| = |\psi_0(e^{i \theta_{n,k}}) - \psi_0(e^{iJ^{-1}(t)})|\tilde{=} |e^{i \theta_{n,k}} - e^{iJ^{-1}(t)}| (|e^{i \theta_{n,k}} + 1| + |e^{i \theta_{n,k}} - e^{iJ^{-1}(t)}|)^{-1/2},
$$
and from ($\ref{ineq:theta0}$), $|\theta_{n,k} - \pi| \ge {{\pi}\over{n+1}}$, it is obvious that $|e^{i \theta_{n,k}} + 1| \tilde{=} |\rho_n e^{i \theta_{n,k}} + 1|$, so that
$$
|\psi_0(e^{i \theta_{n,k}}) - \psi_0(e^{iJ^{-1}(t)})|\tilde{=} |e^{i \theta_{n,k}} - e^{iJ^{-1}(t)}| (|\rho_n e^{i \theta_{n,k}} + 1| + |e^{i \theta_{n,k}} - e^{iJ^{-1}(t)}|)^{-1/2}.
$$
On the other hand, under the condition of  $\theta_{n,k}\in I_4$, it follows from the distance estimation ($\ref{ineq:dist_11}$) that
$$ 
|{{z^*_{n,k} - \zeta^*_{n,k}}}|\tilde{=} {{|\rho_n e^{i \theta_{n,k}} + 1|^{-1/2}}\over{n}},
$$
so that
\begin{eqnarray*} 
& &\sum_{k\neq k_2(t), \theta_{n,k}\in I_4}\big |{{z^*_{n,k} - \zeta^*_{n,k}}\over{z - z^*_{n,k}}}\big|^2\\
&\le &  c_2\sum_{k\neq k_2(t)} {{|\rho_ne^{i \theta_{n,k}} + 1| + |e^{i \theta_{n,k}} - e^{iJ^{-1}(t)}|}\over{n^2|e^{iJ^{-1}(t)} - e^{i \theta_{n,k}}|^2|\rho_n e^{i \theta_{n,k}} + 1|}} \\
& = & c_2\sum_{k\neq k_2(t)} \big({{1}\over{n^2|e^{iJ^{-1}(t)} - e^{i \theta_{n,k}}|^2}} + {{1}\over{n^2|e^{iJ^{-1}(t)} - e^{i \theta_{n,k}}||\rho_n e^{i \theta_{n,k}} + 1|}}\big)\\
&\le&c_2\sum_{k\neq k_2(t)} {{1}\over{n^2|e^{iJ^{-1}(t)} - e^{i \theta_{n,k}}|^2}} + c_2\big(\sum_{k\neq k_2(t)} {{1}\over{n^2|e^{iJ^{-1}(t)} - e^{i \theta_{n,k}}|^2}}\big)^{1/2}\big(\sum_{k\neq k_2(t)} {{1}\over{n^2|\rho_n e^{i \theta_{n,k}} + 1|^2}}\big)^{1/2}.
\end{eqnarray*} 
Combining with ($\ref{ineq:sum_sq00}$) and ($\ref{ineq:sum_sq10}$), we have
\begin{equation}\label{eq:sum_sq0}
\sum_{k\neq k_1(t),k_2(t), \theta_{n,k}\in [0, \pi)}\big |{{z^*_{n,k} - \zeta^*_{n,k}}\over{z - z^*_{n,k}}}\big|^2 
\le  c_2.
\end{equation}

For the sum of $\{k:\theta_{n,k}<0\}$, we can apply the anti-symmetry property of $z^*_{n,k} = \overline{z^*}_{n,2\lfloor{n/2}\rfloor + 1 - k}$. For $\theta_{n,k_1(t)} > 0$, to conclude that
\begin{eqnarray}\label{NegThetaSum00}
\sum_{\theta_{n,k}<0}\big |{{z^*_{n,k} - \zeta^*_{n,k}}\over{z - z^*_{n,k}}}\big|^2 
&=&    \big |{{z^*_{n,k_1(t)} - \zeta^*_{n,k_1(t)}}\over{z - \overline{z^*}_{n,k_1(t)}}}\big|^2  + \big |{{z^*_{n,k_2(t)} - \zeta^*_{n,k_2(t)}}\over{z - \overline{z^*}_{n,k_2(t)}}}\big|^2 + \sum_{k\neq k_1(t),k_2(t), \theta_{n,k}\in [0, \pi)}\big |{{z^*_{n,k} - \zeta^*_{n,k}}\over{z - \overline{z^*}_{n,k}}}\big|^2\nonumber\\
&\le&  \big |{{z^*_{n,k_1(t)} - \zeta^*_{n,k_1(t)}}\over{\psi_0(1) - z^*_{n,k_1(t)}}}\big|^2  + \big |{{z^*_{n,k_2(t)} - \zeta^*_{n,k_2(t)}}\over{\psi_0(-1) - z^*_{n,k_2(t)}}}\big|^2 + c_2.
\end{eqnarray}
From the distance estimation of ($\ref{ineq:dist21}$) and ($\ref{ineq:dist22}$), we also have
\begin{equation}\label{eq:z-znk1}
\sum_{\theta_{n,k}< 0}\big |{{z^*_{n,k} - \zeta^*_{n,k}}\over{z - z^*_{n,k}}}\big|^2 \le c_2, \hskip 15pt\theta_{n,k_1(t)} > 0,
\end{equation}
and for $\theta_{n,k_1(t)} = 0$, similar to $(\ref{NegThetaSum00})$, we have
$$
\sum_{\theta_{n,k}<0}\big |{{z^*_{n,k} - \zeta^*_{n,k}}\over{z - z^*_{n,k}}}\big|^2 
=    \big |{{z^*_{n,k_1(t)} - \zeta^*_{n,k_2(t)}}\over{z - \overline{z^*}_{n,k_2(t)}}}\big|^2 + \sum_{k\neq k_1(t),k_2(t), \theta_{n,k}\in [0, \pi)}\big |{{z^*_{n,k} - \zeta^*_{n,k}}\over{z - \overline{z^*}_{n,k}}}\big|^2 \le c_2.
$$
From ($\ref{eq:sum_sq0}$), we obtain
\begin{equation}\label{eq:sum_sq}
\sum_{k\neq k_1(t),k_2(t)}\big |{{z^*_{n,k} - \zeta^*_{n,k}}\over{z - z^*_{n,k}}}\big|^2 
\le  c_2.
\end{equation}

Now we are ready to estimate the sum of $\big|\operatorname{Re}\big({{z^*_{n,k} - \zeta^*_{n,k}}\over{z - z^*_{n,k}}}\big)\big |$ in the left-hand side of $(\ref{sumLog2Diff1})$. For $k_1(t) > 0$, observe that  
\begin{equation}\label{eq:JustRe0}
\big|\operatorname{Re}\big({{z^*_{n,0} - \zeta^*_{n,0}}\over{z - z^*_{n,0}}}\big)\big| \le c_2, \hskip15pt  k_1(t) \neq 0.
\end{equation}
On the other hand,
\begin{equation}\label{eq:z-znk2}
z^*_{n,k} - \zeta^*_{n,k} =  (1-\rho_n)e^{i\theta_{n,k}}\psi_0'(e^{i\theta_{n,k}}) - 
e^{2i\theta_{n,k}}\int_1^{\rho_n}\int_1^{s}\psi_0''(re^{i\theta_{n,k}})drds.
\end{equation}
Since $\operatorname{Im} (\log(z - \psi_0(e^{i \theta})))$ is a constant for $0\le\theta\le\pi$, $\frac{d(\log(z - \psi_0(e^{i \theta}))}{d\theta}$ is real. This means that
\begin{equation}\label{eq:re-zero}
\operatorname{Re} \big({{e^{i\theta_{n,k}}\psi_0'(e^{i\theta_{n,k}})}\over{z - z^*_{n,k}}}\big) = \operatorname{Im}\big(-\frac{d(\log(z - \psi_0(e^{i \theta}))}{d\theta} \big|_{\theta = \theta_{n,k}}\big ) = 0, \hskip 5pt  0 \le k \le \lfloor{n/2}\rfloor.
\end{equation}
For the integral expression in $(\ref{eq:z-znk2})$, if $k \neq 0$, then
\begin{eqnarray}\label{doubleIntg00}
\int_1^{\rho_n}\int_1^{s}|\psi_0''(re^{i\theta_{n,k}})|drds &\le& {{(1-\rho_n)^2}\over{2}}\max_{1\le r \le \rho_n}|\psi_0''(re^{i\theta_{n,k}})| \nonumber\\
&\tilde{=}& {{1}\over{n^2|e^{i\theta_{n,k}} - 1|^{1/2}}} + {{1}\over{n^2| e^{i\theta_{n,k}} + 1|^{3/2}}}\nonumber\\
&\tilde{=}& {{1}\over{n^2|\rho_n e^{i\theta_{n,k}} - 1|^{1/2}}} + {{1}\over{n^2|\rho_n e^{i\theta_{n,k}} + 1|^{3/2}}}, \hskip15pt k \neq 0. 
\end{eqnarray}
That the last step holds, since the distance of $\theta_{n,k}$ to both $0$ and $\pi$ is at least $\frac{\pi}{n+1}$ for $k>0$. Combining this observation with the inequalities in (\ref{eq:JustRe0}) and (\ref{doubleIntg00}), we have, for $k_1(t) > 0$,
\begin{eqnarray*}
\sum _{k \neq k_1(t),\theta_{n,k}\in I_1}|\operatorname{Re}\big({{z^*_{n,k} - \zeta^*_{n,k}}\over{z - z^*_{n,k}}} \big )|
&\le& |\operatorname{Re}\big({{z^*_{n,0} - \zeta^*_{n,0}}\over{z - z^*_{n,0}}} \big )| + c_2\sum_{k > 0, k \neq k_1(t),\theta_{n,k}\in I_1}  {{1}\over{n^2|\rho_n e^{i\theta_{n,k}} - 1|^{1/2}|z^*_{n,k} - z|}} \\
& \le & c_2\sum_{k \neq k_1(t),\theta_{n,k}\in I_1}  {{1}\over{n^2|\rho_n e^{i\theta_{n,k}} - 1||e^{i\theta_{n,k}} - e^{i t} |}}\\
& \le & c_2\big\{\sum{{1}\over{n^2|\rho_n e^{i\theta_{n,k}} - 1|^2}}\big\}^{1/2}
\big\{\sum_{k \neq k_1(t)}{{1}\over{n^2|e^{i\theta_{n,k}} - e^{i t} |^2}}\big\}^{1/2}, \hskip 5pt k_1(t) > 0.
\end{eqnarray*}
If $k_1(t) = 0$, we also have
\begin{eqnarray*}
\sum _{k \neq k_1(t),\theta_{n,k}\in I_1}|\operatorname{Re}\big({{z^*_{n,k} - \zeta^*_{n,k}}\over{z - z^*_{n,k}}} \big )|
&\le&  c_2\sum_{k > 0, \theta_{n,k}\in I_1}  {{1}\over{n^2|\rho_n e^{i\theta_{n,k}} - 1|^{1/2}|z^*_{n,k} - z|}} \\
& \le & c_2\big\{\sum_{k>0,\theta_{n,k}\in I_1}{{1}\over{n^2|\rho_n e^{i\theta_{n,k}} - 1|^2}}\big\}^{1/2}
\big\{\sum_{k \neq k_1(t)}{{1}\over{n^2|e^{i\theta_{n,k}} - e^{i t} |^2}}\big\}^{1/2}.
\end{eqnarray*}
In summary, we may conclude that
\begin{equation}\label{ineq:re_sum00}
\sum _{k \neq k_1(t),\theta_{n,k}\in I_1}|\operatorname{Re}\big({{z^*_{n,k} - \zeta^*_{n,k}}\over{z - z^*_{n,k}}} \big )| \le c_2.
\end{equation}
For $\theta_{n,k}\in I_2$, we have
$$
|\operatorname{Re}\big({{z^*_{n,k} - \zeta^*_{n,k}}\over{z - z^*_{n,k}}} \big )| = \big |\frac{\int_1^{\rho_n}\int_1^{s}|\psi_0''(re^{i\theta_{n,k}})|drds }{z - z^*_{n,k}}\big | \le  \frac{c_2}{n^2|z - z^*_{n,k}|}.
$$
From the distance estimation in $(\ref{ineq:dist_01})$, we have
\begin{eqnarray}\label{ineq:re_sum01}
\sum _{k \neq k_1(t),\theta_{n,k}\in I_2}|\operatorname{Re}\big({{z^*_{n,k} - \zeta^*_{n,k}}\over{z - z^*_{n,k}}} \big )| &\le& c_2\sum_{k \neq k_1(t),\theta_{n,k}\in I_2}  {{1}\over{n^2|z - z^*_{n,k}|}} \nonumber\\
& \tilde{=} & c_2\sum_{k \neq k_1(t),\theta_{n,k}\in I_2}  {{1}\over{n^2|e^{i\theta_{n,k}} - e^{i t} | (|e^{i\theta_{n,k}} - e^{i {2\over 3} \pi}| +|e^{i\theta_{n,k}} - e^{i t}|)  }}\nonumber\\
& \le & c_2\sum_{k \neq k_1(t)}  {{1}\over{n^2|e^{i\theta_{n,k}} - e^{i t} |^{2}  }}\le  c_2.
\end{eqnarray}
Similarly, we may also conclude that
\begin{equation}\label{ineq:re_sum02}
\sum _{k \neq k_2(t),\theta_{n,k}\in I_3}|\operatorname{Re}\big({{z^*_{n,k} - \zeta^*_{n,k}}\over{z - z^*_{n,k}}} \big )| \le c_2.
\end{equation}
For $\theta_{n,k}\in I_4$ with $k \neq k_2(t)$, we obtain
$$
{{1}\over{|z - z^*_{n,k}|}}\tilde{=} {{(|e^{i\theta_{n,k}} + 1| + |e^{i\theta_{n,k}} - e^{iJ^{-1}(t)}|)^{1/2}}\over {|e^{i\theta_{n,k}} - e^{iJ^{-1}(t)}|}} \tilde{=} {{|e^{i\theta_{n,k}} + 1|^{1/2} + |e^{i\theta_{n,k}} - e^{iJ^{-1}(t)}|^{1/2}}\over {|e^{i\theta_{n,k}} - e^{iJ^{-1}(t)}|}}.
$$ 
Combining the inequalities $(\ref{eq:z-znk2})$ and $(\ref{doubleIntg00})$, we have
\begin{eqnarray}\label{ineq:re_sum03}
\sum _{k \neq k_2(t),\theta_{n,k}\in I_4}|\operatorname{Re}\big({{z^*_{n,k} - \zeta^*_{n,k}}\over{z - z^*_{n,k}}} \big )| &\tilde{=}& \sum_{k \neq k_2(t),\theta_{n,k}\in I_4}  {{1}\over{n^2|\rho_n e^{i\theta_{n,k}} +1|^{3/2}|z - z^*_{n,k}|}}\nonumber \\
& \tilde{=} & \sum_{k \neq k_2(t),\theta_{n,k}\in I_4}  {{|e^{i\theta_{n,k}} + 1|^{1/2} + |e^{i\theta_{n,k}} - e^{iJ^{-1}(t)}|^{1/2}}\over{n^2|\rho_n e^{i\theta_{n,k}} +1|^{3/2}|e^{i\theta_{n,k}} - e^{iJ^{-1}(t)}|}}\nonumber\\
& = & \sum_{k \neq k_2(t),\theta_{n,k}\in I_4}  {{1}\over{n^2|\rho_n e^{i\theta_{n,k}} +1||e^{i\theta_{n,k}} - e^{iJ^{-1}(t)}|}  }\nonumber\\ 
&+& \sum_{k \neq k_2(t),\theta_{n,k}\in I_4}  {{1}\over{n^2|\rho_n e^{i\theta_{n,k}} +1|^{3/2}|e^{i\theta_{n,k}} - e^{iJ^{-1}(t)}|^{1/2}}  }\nonumber\\
&\le&  \big\{\sum_{k=0}^n{{1}\over{n^2|\rho_n e^{i\theta_{n,k}} + 1|^2}}\big\}^{1/2}
\big\{\sum_{k \neq k_2(t)}{{1}\over{n^2|e^{i\theta_{n,k}} - e^{iJ^{-1}(t)} |^2}}\big\}^{1/2}\nonumber\\
&+& \big\{\sum_{k=0}^n{{1}\over{n^2|\rho_n e^{i\theta_{n,k}} + 1|^2}}\big\}^{3/4}
\big\{\sum_{k \neq k_2(t)}{{1}\over{n^2|e^{i\theta_{n,k}} - e^{iJ^{-1}(t)} |^2}}\big\}^{1/4} \le c_2.
\end{eqnarray}

However, we still have to consider the summation for $\{k:-\pi<\theta_{n,k}<0\}$; that is, $\lfloor{n/2}\rfloor < k \le n$,  but do not have the equivalence of ($\ref{eq:re-zero}$) to rely on. So, let us consider the following argument. Since $0 = \psi_0(1)\in \CC$ and both $z^*_{n,k}$ and $0$ lie on the same line segment of the $L$-shape open arc $\gamma_0$, similar to (\ref{eq:re-zero}) 
$$
\operatorname{Re}\big({{ e^{i \theta_{n,k}}\psi_0'(e^{i \theta_{n,k}})}\over{(\psi_0(1) - z^*_{n,k})}} \big) = 0, \hskip 10pt \lfloor{n/2}\rfloor < k \le n,
$$
so that
\begin{eqnarray*}
\big | \operatorname{Re}{{ e^{i \theta_{n,k}}\psi_0'(e^{i \theta_{n,k}})}\over{(z - z^*_{n,k})}} \big| &=& 
\big | \operatorname{Re}\big ( {{ e^{i \theta_{n,k}}\psi_0'(e^{i \theta_{n,k}})}\over{(\psi_0(1) - z^*_{n,k})}} {{{\psi_0(1) - z^*_{n,k}}\over{z - z^*_{n,k}}}}\big )  \big| \\
&=&\big | {{ e^{i \theta_{n,k}}\psi_0'(e^{i \theta_{n,k}})}\over{(\psi_0(1) - z^*_{n,k})}} \operatorname{Im}{{\psi_0(1) - z^*_{n,k}}\over{z - z^*_{n,k}}}  \big| \\
&=&\big | {{ e^{i \theta_{n,k}}\psi_0'(e^{i \theta_{n,k}})}\over{(\psi_0(1) - z^*_{n,k})}} \operatorname{Im}\big ( {{\psi_0(1) - z}\over{z - z^*_{n,k}}}+1\big )  \big|\\
& = & {{|\psi_0'(e^{i \theta_{n,k}})|}\over {|z^*_{n,k}|}}{{{|\operatorname{Im}(-z(\overline{z} - \overline{ z^*_{n,k}})|}\over {|z - z^*_{n,k}|^2}}} \\
 & = & {{|\psi_0'(e^{i \theta_{n,k}})| |z|}\over {|z|^2 + |z^*_{n,k}|^2}}, \hskip 10pt \lfloor{n/2}\rfloor < k \le n.
\end{eqnarray*}
Then it follows that
\begin{equation}
\sum_{k = \lfloor{n/2}\rfloor + 1}^n \big | \operatorname{Re}{{ (\rho_n - 1)e^{i \theta_{n,k}}\psi_0'(e^{i \theta_{n,k}})}\over{(z - z^*_{n,k})}} \big| = {{1}\over{n+1}} \sum_{\theta_{n,k}<0}{{|\psi_0'(e^{i \theta_{n,k}})| |z|}\over {|z|^2 + |\psi_0(e^{i \theta_{n,k}})|^2}}.
\end{equation}
For any $s\in [\theta_{n,k}, \theta_{n,k+1}]$, we have
$$
{{|\psi_0'(e^{i \theta_{n,k}})| |z|}\over {|z|^2 + |\psi_0(e^{i \theta_{n,k}})|^2}}\tilde{=}{{|\psi_0'(e^{i s})| |z|}\over {|z|^2 + |\psi_0(e^{i s})|^2}},
$$
so that
\begin{eqnarray*}
\sum_{k = \lfloor{n/2}\rfloor + 1}^n \big | \operatorname{Re}{{ (\rho_n - 1)e^{i \theta_{n,k}}\psi_0'(e^{i \theta_{n,k}})}\over{(z - z^*_{n,k})}} \big| & \tilde{=} & \frac{1}{2\pi}\int_{\theta_{n,k}}^{\theta_{n,k+1}} {{|\psi_0'(e^{i s})| |z|}\over {|z|^2 + |\psi_0(e^{i s})|^2}} ds\\
& = & \frac{1}{2\pi}\int_0^{\pi} {{|z|}\over {|z|^2 + |\psi_0(e^{i s})|^2}} |d\psi_0(e^{i s})| \\
&=& {{1}\over{\pi}} \arctan({{\sqrt[4]{27}}\over{|z|}})\\
&\le & {1\over2}.
\end{eqnarray*}
From ($\ref{eq:z-znk2}$),
\begin{equation}\label{ineq:re_sum2}
\sum_{k = \lfloor{n/2}\rfloor + 1}^n|\operatorname{Re}\big({{z^*_{n,k} - \zeta^*_{n,k}}\over{z - z^*_{n,k}}} \big )| \le \sum_{k = \lfloor{n/2}\rfloor + 1}^n \big(\big | \operatorname{Re}{{ (\rho_n - 1)e^{i \theta_{n,k}}\psi_0'(e^{i \theta_{n,k}})}\over{(z - z^*_{n,k})}} \big| + {{\int_1^{\rho_n}\int_1^{s}|\psi_0''(re^{i\theta_{n,k}})|drds}\over{|\psi_0(1) - z^*_{n,k}|}}\big)\le c_2.
\end{equation}

\noindent Combining with ($\ref{eq:sep_sum}$), ($\ref{eq:sum_sq}$), ($\ref{ineq:re_sum01}$), ($\ref{ineq:re_sum02}$), ($\ref{ineq:re_sum03}$) and ($\ref{ineq:re_sum2}$), we have
\begin{equation}
\sum_{k\neq k_1(t), k_2(t)}\big |\big |{{z - \zeta^*_{n,k}}\over{z - z^*_{n,k}}}\big|^2 - 1\big |
\le   c_2.
\end{equation}
This implies that ($\ref{ineq:Omega_ToOmega1}$) holds.

That ($\ref{ineq:Key0}$) holds is a consequence of ($\ref{ineq:Omega_ToOmega1}$) and ($\ref{ineq:omega_bound2}$) in Lemma $\ref{lemma:omega_bound2}$. From ($\ref{ineq:Key0}$), we have
\begin{equation}
\big |{{\omega_n(z)}\over{z - z^*_{n,k_1(t)}}}\big | \tilde{=} \big |{{  z - z^*_{n,k_2(t)} } \over {(z - \zeta^*_{n,k_1(t)}) ( z - \zeta^*_{n,k_2(t)})}}\big|, \hskip 15pt z = \psi_0(e^{i t}), 0\le |t|\le{{2\pi}\over3}.
\end{equation}
The proof of the lemma is complete by considering the limit of $z\rightarrow z^*_{n,j}$, or equivalently $t\rightarrow \theta_{n,j}$, yielding ($\ref{ineq:Key'01}$).

\end{Proof}
\vskip .1in
For the family ${\bf{z^{**}_n}} = \{z^{**}_{n,k}\}_{k =0}^n$ of adjusted Fej\'er points, we set
\begin{equation} \label{ineq:Omega-tilde}
  \tilde\omega_n(z) := \prod_{0\le k \le n} (z - z^{**}_{n,k}),
\end{equation}
and must show that the family of points $\{\zeta^*_{n,k}\}_{k =0}^n$ on the level curve $\Gamma_n$ can be used to replace ${\bf{z^{**}_n}} = \{z^{**}_{n,k}\}_{k =0}^n$ on the curve $\gamma_0$. For this purpose, we will first show that $\tilde\omega_n(z)$ has similar properties as those of $\omega_n(z)$ in Lemma $\ref{lemma:Omega_ToOmega1}$. 

\begin{lemma}\label{lemma:Key00}
Let  ${\bf{z}^{**}_n} = \{z^{**}_{n,k}\}_{k =0}^n$ be the adjustment of the Fej\'er points $\{z^{*}_{n,k}\}_{k =0}^n$, that lie on the open arc $\gamma_0$. Then for $z = \psi_0(e^{it})\in\gamma_0$ with $ 0\le |t|\le{{2\pi}\over3} $,
\begin{equation}\label{equivalence}
\begin{cases}
{\displaystyle |\tilde\omega_n(z)| \tilde{=} \big |{{ (z - z^{**}_{n,k_1(t)}) ( z - z^{**}_{n,k_2(t)}) } \over {(z - \zeta^*_{n,k_1(t)}) ( z - \zeta^*_{n,k_2(t)})}}\big|};\\
&\\  {\displaystyle |\tilde\omega'_n(z^{**}_{n,j})| \tilde{=} \big |{{ z^{**}_{n,j} - z^{**}_{n,k_2(\tilde\theta^*_{n,j})} } \over {(z^{**}_{n,j} - \zeta^*_{n,j}) ( z^{**}_{n,j} - \zeta^*_{n,k_2(\tilde\theta^*_{n,j})}})}\big|   , \hskip15pt |\theta_{n,j}| \le{{2\pi}\over3}}; \\
&\\  {\displaystyle |\tilde\omega'_n(z^{**}_{n,j})| \tilde{=} \big |{{ z^{**}_{n,j} - z^{**}_{n,k_1(\tilde\theta^*_{n,j})} } \over {(z^{**}_{n,j} - \zeta^*_{n,j}) ( z^{**}_{n,j} - \zeta^*_{n,k_1(\tilde\theta^*_{n,j})}})}\big|\tilde{=}{{1} \over {|z^{**}_{n,j} - \zeta^*_{n,j}|}}, \hskip 15pt |\theta_{n,j}| >{{2\pi}\over3}.}\\
\end{cases}
\end{equation}
\end{lemma}

\begin{Proof}
Let us first prove that there exist constants $0<c_1\le c_2$, such that for all $z= \psi_0(e^{it})\in \gamma_0$ with $0\le t\le\frac{2\pi}{3}$,
\begin{equation} \label{ineq:omega_ToOmega2}
c_1\le \prod_{j\neq k_1(t), k_2(t)}{{|z - z^*_{n,j}|}\over {|z - z^{**}_{n,j}|}}\le c_2.
\end{equation}
To prove (\ref{ineq:omega_ToOmega2}), observe that for $j\neq k_1(t)$, since $|t - \theta_{n,j}| \ge \pi/(n+1)$, in view of (\ref{adjDiffJK}), it is clear that 
$$
|t - \tilde\theta^*_{n.j}| \le |t - \theta_{n.j}| + |\theta_{n,j} - \tilde\theta^*_{n.j}| \le  (1 + 2/3) |t - \theta_{n,j}|.
$$
On the other hand, it is also clear that
\begin{equation}\label{t_theta*0}
(1 - 2/3) |t - \theta_{n,j}| \le |t - \tilde\theta^*_{n.j}| \le (1 + 2/3) |t - \theta_{n,j}|, \hskip 15pt j\neq k_1(t), |\theta_{n,j}|\le {2\pi}/3.
\end{equation}
Similarly, for $j\neq k_2(t)$, by ($\ref{adjDiffJK}$), we also have 
\begin{equation}\label{t_theta_1}
(1 - 2/3) |J^{-1}(t) - \theta_{n,j}| \le |J^{-1}(t) - \tilde\theta^*_{n.j}| \le (1 + 2/3) |J^{-1}(t) - \theta_{n,j}|, \hskip 15pt j\neq k_2(t), |\theta_{n,j}|\ge {2\pi}/3.
\end{equation}
From ($\ref{ineq:dist_00}$), ($\ref{ineq:dist_01}$), ($\ref{ineq:dist_10}$), ($\ref{ineq:dist_11}$) and ($\ref{ineq:dist_12}$), regardless the cases, we always have
\begin{equation}\label{t_theta_2}
|z - z^{**}_{n,j}| \tilde{=} |z - z^*_{n,j}|, \hskip 15pt j\neq k_1(t), j\neq k_2(t).
\end{equation}
On the other hand, it is obvious that
$$
\prod_{j = 0}^n{{|z - z^*_{n,j}|}\over {|z - z^{**}_{n,j}|}} = \prod_{\{k:|\theta_{n,k}| > \frac{2\pi}{3}, (j(k),k)_n \text{ adjustment pair }\}}\big|{{(z - z^*_{n,k})(z - z^*_{n,j(k)})}\over{(z - z^{**}_{n,k})(z - z^{**}_{n,j(k)})}}\big |.
$$
From $(\ref{t_theta_2})$, we may ignore several finite terms that do not affect the estimation result of ($\ref{ineq:omega_ToOmega2}$) to conclude that
\begin{equation}\label{ineq:pair00}
\prod_{j\neq k_1(t), k_2(t)}{{|z - z^*_{n,j}|}\over {|z - z^{**}_{n,j}|}} \tilde{=} \prod_{k \in \M_n(t)} \big |{{(z - z^*_{n,k})(z - z^*_{n,j(k)})}\over{(z - z^{**}_{n,k})(z - z^{**}_{n,j(k)})}}\big |,
\end{equation}
where $\M_n(t) := \{k:|\theta_{n,k}| > \frac{2\pi}{3}, (j(k),k)_n \text{ adjustment pair, } k\neq k_2(t), j(k) \neq k_1(t)\}$.

Next from the identity:
\begin{equation}\label{ineq:sepTo1+}
{{(z - z^*_{n,k})(z - z^*_{n,j(k)})}\over{(z - z^{**}_{n,k})(z - z^{**}_{n,j(k)})}} = 1 + { {z^{**}_{n,k}+ z^{**}_{n,j(k)} - z^*_{n,k} - z^*_{n,j(k)}} \over{z - z^{**}_{n,j(k)}} } + {  {(z^*_{n,k} - z^{**}_{n,k})(z^*_{n,j(k)} - z^{**}_{n,k}) }\over{(z - z^{**}_{n,k})(z - z^{**}_{n,j(k)})}},
\end{equation}
we obtain 
\begin{eqnarray*}
& &\big |{ {z^{**}_{n,k}+ z^{**}_{n,j(k)} - z^*_{n,k} - z^*_{n,j(k)}} }\big | \\
&=& \big | \psi_0(e^{i\tilde\theta^*_{n,k}}) - \psi_0(e^{i J(\theta_{n,k})}) + \psi_0(e^{i\tilde\theta^*_{n,j(k)}}) - \psi_0(e^{i \theta_{n,j(k)}}) \big | \\
& = & \big |{\int_{J(\theta_{n,k})}^{\tilde\theta^*_{n,k}} e^{is} \psi_0'(e^{is}) ds  +  \int_{\theta_{n,j(k)}}^{\tilde\theta^*_{n,j(k)}} e^{is} \psi_0'(e^{is}) ds} \big |\\
& = & \big|\int_{0}^{\tilde\theta^*_{n,k} - J(\theta_{n,k})} e^{i(s+J(\theta_{n,k}))}\psi_0'(e^{i(s+J(\theta_{n,k}))}) ds
- \int_{0}^{\theta_{n,j(k)} - \tilde\theta^*_{n,j(k)}} e^{i(s + \tilde\theta^*_{n,j(k)})}\psi_0'(e^{i(s + \tilde\theta^*_{n,j(k)})}) ds \big|.
\end{eqnarray*}
Also, from the definition $(\ref{ineq:thetaTildeJ})$ and $(\ref{ineq:thetaTildeK})$, we have
\begin{equation}\label{ineq:tilde_theta1}
\theta_{n,j(k)} + J(\theta_{n,k}) = \tilde\theta^*_{n,j(k)} + \tilde\theta^*_{n,k}. 
\end{equation}
This implies that ${\tilde\theta^*_{n,k} - J(\theta_{n,k})} = \theta_{n,j(k)} - \tilde\theta^*_{n,j(k)}$; and hence,
\begin{eqnarray*}
& &\big |{ {z^{**}_{n,k}+ z^{**}_{n,j(k)} - z^*_{n,k} - z^*_{n,j(k)}}  }\big |\\
&=&\big |{{\int_{0}^{\tilde\theta^*_{n,k} - J(\theta_{n,k})} \int_{s+\tilde\theta^*_{n,j(k)}}^{s+J(\theta_{n,k})} \big(e^{i\eta} \psi_0'(e^{i\eta})\big)' d\eta ds }} \big |\\
&\le & {{|{\tilde\theta^*_{n,k} - J(\theta_{n,k}})| |J(\theta_{n,k}})-{\tilde\theta^*_{n,j(k)} | \max_{\eta \in [\tilde\theta^*_{n,k}, \tilde\theta^*_{n,j(k)}]} \big(|\psi_0'(e^{i\eta})|  + |\psi_0''(e^{i\eta})|\big) }}\\
&\le & c_2 {{|\psi_0''(e^{i\tilde\theta^*_{n,j(k)}})|} \over{n^2}},
\end{eqnarray*}
where the last step holds, because of $|\tilde\theta^*_{n,j(k)}| \tilde{=} |\tilde\theta^*_{n,k}|$,  for $|\theta_{n,j(k)} - J(\theta_{n,k})| \le \pi/(n+1)$.

From the distance estimation, for $j(k) \neq k_1(t)$ and $|t - \tilde\theta^*_{n,j(k)}| \ge \frac{\pi}{n+1}$, we obtain
$$
|z - z^{**}_{n,j(k)}|  \tilde{=} |z - \psi_0(\rho_n e^{i\tilde\theta^*_{n,j(k)}})| \ge c_1|\psi_0'(\rho_n e^{i\tilde\theta^*_{n,j(k)}})||t - \tilde\theta^*_{n,j(k)}|. 
$$
This implies that
\begin{eqnarray*}
\big |{ {z^{**}_{n,k}+ z^{**}_{n,j(k)} - z^*_{n,k} - z^*_{n,j(k)} \over {z - z^{**}_{n,j(k)}}}  }\big | 
&\le & c_2 {{|\psi_0''(e^{i\tilde\theta^*_{n,j(k)}})|} \over{n^2|\psi_0'(\rho_n e^{i\tilde\theta^*_{n,j(k)}})||t - \tilde\theta^*_{n,j(k)}|}}.
\end{eqnarray*}
Also, for $|\tilde\theta^*_{n,j(k)}| \le {{2\pi}\over{3}}$, since $|\psi_0''( e^{i\tilde\theta^*_{n,j(k)}})| \tilde{=} |\tilde\theta^*_{n,j(k)}|^{-1/2}$ and
$$|\psi_0'(\rho_n e^{i\tilde\theta^*_{n,j(k)}})| \tilde{=}  ({{2\pi}\over{3}} - |\tilde\theta^*_{n,j(k)}|+ \frac1{n})|\tilde\theta^*_{n,j(k)} + \frac1{n}|^{1/2},$$  we have
\begin{eqnarray}\label{sum001}
& &\sum_{k\in \M_(t)}\big |{ {z^{**}_{n,k}+ z^{**}_{n,j(k)} - z^*_{n,k} - z^*_{n,j(k)}} \over{z - z^{**}_{n,j(k)}} }\big | \nonumber\\
&\le & c_2 \sum_{j(k) \neq k_1(t)}{{1} \over{n^2|t - \tilde\theta^*_{n,j(k)}|({{2\pi}\over{3}} - |\tilde\theta^*_{n,j(k)}|+\frac1{n})|\tilde\theta^*_{n,j(k)}|}} \nonumber\\
&\le & c_2 \big(\sum_{j(k) \neq k_1(t)}{{1} \over{n^2|t - \tilde\theta^*_{n,j(k)}|^2}}\big)^{1/2}\big(\sum{{1} \over{n^2({{2\pi}\over{3}} - |\tilde\theta^*_{n,j(k)}|+\frac1{n})^2|\tilde\theta^*_{n,j(k)}|^2}}\big)^{1/2} \nonumber\\
&\le & c_2 \big(\sum_{|\tilde\theta^*_{n,j(k)}| > {{2\pi}\over{3}}} {{1} \over{n^2({{2\pi}\over{3}} - |\tilde\theta^*_{n,j(k)}|+\frac1{n})^2}} + \sum_{|\tilde\theta^*_{n,j(k)}| \le {{2\pi}\over{3}}}{{1} \over{n^2|\tilde\theta^*_{n,j(k)}|^2}}\big)^{1/2} \nonumber\\
&\le & c_2\big(1 + \sum{{1} \over{n^2|\tilde\theta^*_{n,j(k)}|^2}}\big)^{1/2} 
\le  c_2,
\end{eqnarray}
where the last step holds, because of $|\tilde\theta^*_{n,j(k)}| \ge J(\theta_{n,\lfloor{n/2}\rfloor}) - \frac{2\pi}{3(n+1)} \tilde{=} \frac1{\sqrt[3]{n}}$ from (\ref{ineq:J(t)ClosePi}). 
\vskip .1in
On the other hand, for $j(k) \neq k_1(t)$ and $k \neq k_2(t)$, since
$$
   \big|{  {(z^*_{n,k} - z^{**}_{n,k})(z^*_{n,j(k)} - z^{**}_{n,k}) }\over{(z - z^{**}_{n,k})(z - z^{**}_{n,j(k)})}} \big | \le c_2 \big| {{\psi_0'(e^{i \tilde\theta^*_{n,j(k)}})}\over {n}} \big | ^2 \big| {{1}\over {(e^{it} - e^{i\tilde\theta^*_{n,j(k)}})\psi_0'(\rho_n e^{i\tilde\theta^*_{n,j(k)}})}} \big | ^2 =  {{c_2}\over {n^2|e^{it} - e^{i\tilde\theta^*_{n,j(k)}}|^2}},
$$
it follows that
\begin{equation}\label{sum002}
    \sum_{k\in \M_n(t)}\big|{  {(z^*_{n,k} - z^{**}_{n,k})(z^*_{n,j(k)} - z^{**}_{n,k}) }\over{(z - z^{**}_{n,k})(z - z^{**}_{n,j(k)})}} \big | \le c_2.\\
\end{equation}
Therefore, by applying ($\ref{ineq:pair00}$), ($\ref{ineq:sepTo1+}$), ($\ref{sum001}$) and ($\ref{sum002}$), we have completed the proof of ($\ref{ineq:omega_ToOmega2}$).\\

To derive the first item of (\ref{equivalence}), we simply apply ($\ref{ineq:omega_ToOmega2}$) and ($\ref{ineq:Key0}$) in Lemma $\ref{lemma:Omega_ToOmega1}$ to conclude that
$$
\big|{{\tilde\omega_n(z)}\over{\omega_n(z)}}\big| = \prod_{0\le j \le n}{{|z - z^*_{n,j}|}\over {|z - z^{**}_{n,j}|}} \tilde{=} \big |{{ (z - z^{**}_{n,k_1(t)}) ( z - z^{**}_{n,k_2(t)}) } \over {(z - z^*_{n,k_1(t)}) ( z - z^*_{n,k_2(t)})}}\big|.
$$
The proofs of the second and third items of equivalence relations in (\ref{equivalence}) are similar to that of ($\ref{ineq:Key'01}$) in Lemma \ref{lemma:Omega_ToOmega1}. This completes the proof of the lemma.
\end{Proof}

\begin{lemma} \label{lemma:localSum}
Let $\{z^{**}_{n,k}\}_{k = 0}^n$ be a family of adjusted Fej\'er points that lie on $\gamma_0$. Then there exists some constant $c_2 >0$ such that
\begin{equation} \label{ineq:SumOmegaOverOmegaPrime01}
\sum_{\{j:0\le |\theta_{n,j}| \le {{2\pi}\over3},  |k_2(\theta_{n,j}) - k_2(t)| \le 1\}} \big|{{\tilde\omega_n(z)}\over{{\tilde\omega'_n(z^{**}_{n,j}) (z - z^{**}_{n,j}) }   }}\big|\le c_2 \log(n),
\end{equation}
for all $z = \psi_0(e^{it})\in\gamma_0$ with $0\le |t|\le{{2\pi}\over3}$. 
\end{lemma}
\begin{Proof}
For any ${2\over 3}\pi \le \theta \le \pi$, it follows from (1.1) that
$$
|\psi_0'(\rho_n e^{i \theta})|\tilde{=} |-1 - \rho_n e^{i \theta}|^{1\over 2} |e^{i {2\over 3} \pi} - \rho_n e^{i \theta}|.
$$
Since $|k_2(\tilde\theta^*_{n,j}) - k_2(t)| \le 1$, we have
$$
|\psi_0'(\rho_n e^{i \theta_{n,k_2(\tilde\theta^*_{n,j})}})| \tilde{=} |\psi_0'(\rho_n e^{i \theta_{n,k_2(t)}})|. 
$$
In addition, it follows from the distance estimation ($\ref{ineq:dist21}$) that
$$
|z - \zeta^*_{n,k_2(t)}| \tilde{=} |e^{i J^{-1}(t)} - \rho_n e^{i \theta_{n,k_2(t)}}| |\psi_0'(\rho_n e^{i \theta_{n,k_2(t)}})|\tilde{=}{{|\psi_0'(\rho_n e^{i \theta_{n,k_2(t)}})|}\over{n}} \tilde{=} {{|\psi_0'(\rho_n e^{i \theta_{n,k_2(\tilde\theta^*_{n,j})}})|}\over{n}} \tilde{=}|z^{**}_{n,j} - \zeta^*_{n,k_2(\tilde\theta^*_{n,j})}|. 
$$ 
By applying ($\ref{equivalence}$) in Lemma $\ref{lemma:Key00}$, we obtain
\begin{eqnarray*}
\big|{{\tilde\omega_n(z)}\over{{\tilde\omega'_n(z^{**}_{n,j}) (z - z^{**}_{n,j}) }   }}\big| &\tilde{=}&
\big|{{ (z - z^{**}_{n,k_1(t)}) ( z - z^{**}_{n,k_2(t)})(z^{**}_{n,j} - \zeta^*_{n,j}) ( z^{**}_{n,j} - \zeta^*_{n,k_2(\tilde\theta^*_{n,j})}) } \over {(z - \zeta^*_{n,k_1(t)}) ( z - \zeta^*_{n,k_2(t)})(z^{**}_{n,j} - z^{**}_{n,k_2(\tilde\theta^*_{n,j})})(z - z^{**}_{n,j})}}\big|\\
&\tilde{=}&\big|{{ (z - z^{**}_{n,k_1(t)}) ( z - z^{**}_{n,k_2(t)})(z^{**}_{n,j} - \zeta^*_{n,j}) } \over {(z - \zeta^*_{n,k_1(t)}) (z^{**}_{n,j} - z^{**}_{n,k_2(\tilde\theta^*_{n,j})})(z - z^{**}_{n,j})}}\big|\\
&=&\big|{{ z - z^{**}_{n,k_1(t)} } \over {z - \zeta^*_{n,k_1(t)}}}\big|\big|{{ (z - z^{**}_{n,k_2(t)})( z^{**}_{n,j} - \zeta^*_{n,j}) } \over { (z^{**}_{n,j} - z^{**}_{n,k_2(\tilde\theta^*_{n,j})})(z - z^{**}_{n,j})}}\big| \\
&\le&c_2\big|{{ z - z^{**}_{n,k_1(t)} } \over {z - \zeta^*_{n,k_1(t)}}}\big| \big|{{  (z - z^{**}_{n,k_2(\tilde\theta^*_{n,j})})( z^{**}_{n,j} - \zeta^*_{n,j}) } \over { (z^{**}_{n,j} - z^{**}_{n,k_2(\tilde\theta^*_{n,j})})(z - z^{**}_{n,j})}}\big|\\
&=&c_2\big|{{ z - z^{**}_{n,k_1(t)} } \over {z - \zeta^*_{n,k_1(t)}}}\big| \big|{{  z^{**}_{n,j} - \zeta^*_{n,j}} \over { z^{**}_{n,j} - z^{**}_{n,k_2(\tilde\theta^*_{n,j})}}} + {{  z^{**}_{n,j} - \zeta^*_{n,j} } \over { z - z^{**}_{n,j}}}\big|,
\end{eqnarray*}
and according to the distance estimation ($\ref{ineq:dist21}$), we have
\begin{eqnarray*}
& &\sum_{\{j:0\le |\theta_{n,j}| \le {{2\pi}\over3},  |k_2(\theta_{n,j}) - k_2(t)| \le 1\}} \big|{{\tilde\omega_n(z)}\over{{\tilde\omega'_n(z^{**}_{n,j}) (z - z^{**}_{n,j}) }   }}\big|\\
&\le&c_2\sum_{\{j:|k_2(\tilde\theta^*_{n,j}) - k_2(t)| \le 1 \}} \big|{{ z^{**}_{n,j} - \zeta^*_{n,j}} \over { z^{**}_{n,j} - z^{**}_{n,k_2(\tilde\theta^*_{n,j})}}}\big| + c_2\sum_{j \neq k_1(t)}\big|{{  z^{**}_{n,j} - \zeta^*_{n,j} } \over { z - z^{**}_{n,j}}}\big|
+ c_2 \big|{{ z - z^{**}_{n,k_1(t)} } \over {z - \zeta^*_{n,k_1(t)}}}\big|\big|{{  z^{**}_{n,k_1(t)} - \zeta^*_{n,k_1(t)} } \over { z - z^{**}_{n,k_1(t)}}}\big|\\ 
&\le& c_2 \sum_{\{j:|k_2(\tilde\theta^*_{n,j}) - k_2(t)| \le 1 \}} {{1} \over { n|\tilde\theta^*_{n,j} - J(\tilde\theta^*_{n,k_2(\tilde\theta^*_{n,j})})|}} + c_2\sum_{j \neq k_1(t)} {{1} \over { n|t - \tilde\theta^*_{n,j}|}} + c_2 \big|{{ z^{**}_{n,k_1(t)} - \zeta^*_{n,k_1(t)} } \over {z - \zeta^*_{n,k_1(t)}}}\big|\\
&\le& c_2\sum_{m = -1}^1\sum_{j = 0}^{n}{{1} \over { n|\tilde\theta^*_{n,j} - J(\tilde\theta^*_{n,k_2(t)+m})|}} + c_2\sum_{j \neq k_1(t)} {{1} \over { n|t - \tilde\theta^*_{n,j}|}} + c_2.
\end{eqnarray*}
This implies that ($\ref{ineq:SumOmegaOverOmegaPrime01}$) and completes the proof of the lemma.
\end{Proof}
\vskip.1in
\bhag{ Proofs of main theorems}
With the preparations in the previous two sections, we are now ready to prove  Theorems 1--3 as stated in the introduction section.
\vskip.1in
\noindent Let us first establish {\bf Theorem 1}.
\vskip.1in
\begin{Proof}  

Let $t_0 = (\theta_{n,0}+\theta_{n,1})/2, z_0 =\psi_0(e^{it_0}) \in\gamma_0$. From $(\ref{k1k2})$, it follows that $k_1(t_0) = 0$ and $k_2(t_0) = \lfloor{n/2}\rfloor$; and from ($\ref{ineq:theta0}$), we have
$$
\theta_{n,k_2(t_0)} = \pi - {{\pi}\over {n+1}}. 
$$
Also, in view of ($\ref{ineq:J(t)ClosePi}$) we obtain 
$J(\theta_{n,k_2(t_0)}) = \sqrt[3]{{4\pi}\over {n+1}} + O(1/n)$. 
Similarly, we have
\begin{equation}\label{ineq:theta^*00}
J(\theta_{n,k}) \tilde{=} \sqrt[3]{{{n + 1 - 2k}\over{n+1}}}, \hskip 15.0pt {{2\pi}\over{3}}<\theta_{n,k}<\pi.
\end{equation} 
From  ($\ref{ineq:Key0}$) in Lemma $\ref{lemma:Omega_ToOmega1}$ and the distance estimation, it follows that
\begin{equation}\label{ineq:BoundOmegaZ0}
|{\omega}_n(z_0)|\tilde{=} \big|{{ (z_0 - z^*_{n,0}) ( z - z^*_{n,\lfloor{n/2}\rfloor}) } \over {(z - \zeta^*_{n,0}) ( z - \zeta^*_{n,k_2(t)})}}\big|\tilde{=} \big | \frac{z - z^*_{n,\lfloor{n/2}\rfloor}}{z - \zeta^*_{n,k_2(t)})}\big |\tilde{=} 1.
\end{equation} 
Therefore,
\begin{equation}\label{ineq:L_nOmega0}
 \sum_{j = 0}^n \big | {{{\omega}_n(z_0)} \over{{\omega}'_n(z^*_{n,j})(z_0 - z^*_{n,j}) }}\big | \ge c_1 \sum_{j = 0}^{n/6}  {{|z^*_{n,j} - \zeta^*_{n,k_2(\theta_{n,j})}||z^*_{n,j} - \zeta^*_{n,j}|}\over{|(z_0-z^*_{n,j})(z^*_{n,j} - z^*_{n,k_2(\theta_{n,j})}})|}. 
\end{equation} 
From the distance estimation $(\ref{ineq:dist_11})$, we have
\begin{eqnarray*}
| z^*_{n,j} - \zeta^*_{n,k_2(\theta_{n,j})} | &=& |\psi_0(e^{i J^{-1}(\theta_{n,j})}) - \psi_0(\rho_n e^{i \theta_{n,k_2(\theta_{n,j})}}) | \\
&\tilde{=}& {|\psi_0'(\rho_n e^{i \theta_{n,k_2(\theta_{n,j})}})|\over n} \\
&\tilde{=}& {1\over {n| 1 + e^{i  \theta_{n,k_2(\theta_{n,j})} } |^{1/2} }}.
\end{eqnarray*}
Also from $(\ref{ineq:dist_00})$, we obtain
$$
| z^*_{n,j} - \zeta^*_{n,j} | = |\psi_0(e^{i \theta_{n,j}}) - \psi_0(\rho_n e^{i \theta_{n,j}}) | \tilde{=} {{|\psi_0'(\rho_n e^{i \theta_{n,j}})|}\over {n}} \tilde{=}{{(1 + j)^{1/2}}\over {n^{3/2}}};
$$
$$
| z^*_{n,j} - z_0 | = |\psi_0(e^{i \theta_{n,j}}) - \psi_0(e^{i t_0}) |\tilde{=} |\theta_{n,j} - t_0| \sqrt{\theta_{n,j} + t_0}\tilde{=} {{(1 + j)^{3/2}}\over {n^{3/2}}};
$$
$$
| z^*_{n,j} - z^*_{n,k_2(\theta_{n,j})} | \tilde{=} |e^{i\theta_{n,j}} - e^{i J(\theta_{n,k_2(\theta_{n,j})})}| |1 - e^{i J(\theta_{n,k_2(\theta_{n,j})})}|^{1/2}. 
$$
Hence, we may conclude that
\begin{equation}\label{ineq:L_nOmega1}
 {{|z^*_{n,j} - \zeta^*_{n,k_2(\theta_{n,j})}||z^*_{n,j} - \zeta^*_{n,j}|}\over{|(z_0-z^*_{n,j})(z^*_{n,j} - z^*_{n,k_2(\theta_{n,j})}})|} \ge   {{c_1}\over{n(1+j)|e^{i\theta_{n,j}} - e^{i J(\theta_{n,k_2(\theta_{n,j})})}| |1 - e^{i J(\theta_{n,k_2(\theta_{n,j})})}|^{1/2}| 1 + e^{i  \theta_{n,k_2(\theta_{n,j})} } |^{1/2}}}. 
 \end{equation} 
In addition, by defining  
\begin{equation} \label{S_K}
S_n(k) := \{j: k_2(\theta_{n,j}) = k\},  \hskip 10pt |\theta_{n,k}| \ge \frac{2\pi}{3}, 
\end{equation}
it follows that
$$(j+1) \le {{(n+1) J(\theta_{n,k - 1})} \over{2\pi}}\tilde{=} {{(n+1) J(\theta_{n,k})}}, \hskip 15pt j\in S_n(k), j\le n/6,
$$
and from  (\ref{ineq:L_nOmega0}), we have
\begin{eqnarray}\label{ineq:L_nOmega2}
 & &\sum_{j = 0}^n \big | {{{\omega}_n(z_0)} \over{{\omega}'_n(z^*_{n,j})(z_0 - z^*_{n,j}) }}\big |\nonumber\\
 & \ge &
 \sum_{k ={\lfloor{5n/12}\rfloor}}^{{\lfloor{n/2}\rfloor}} \hskip 10pt\sum_{j\in S_n(k)} {{|z^*_{n,j} - \zeta^*_{n,k_2(\theta_{n,j})}||z^*_{n,j} - \zeta^*_{n,j}|}\over{|(z_0-z^*_{n,j})(z^*_{n,j} - z^*_{n,k_2(\theta_{n,j})}})|} \nonumber\\ 
 &\ge& c_1 \sum_{k ={\lfloor{5n/12}\rfloor}}^{{\lfloor{n/2}\rfloor}}{{1}\over{nJ(\theta_{n,k})|1 - e^{i J(\theta_{n,k})}|^{1/2}| 1 + e^{i  \theta_{n,k} } |^{1/2}}}\sum_{j\in S_n(k)} {{1}\over{n|e^{i\theta_{n,j}} - e^{i J(\theta_{n,k})}| }}. 
\end{eqnarray}
Also, it follows from ($\ref{ineq:J'(t)ClosePi}$) that 
$$
 |J(\theta_{n,k}) -  J(\theta_{n,k-1})| \tilde{=} {{2\pi|J'(\theta_{n,k})|}\over{n+1}}\tilde{=} {{2\pi}\over{n+1}}\big({{n+1}\over{n + 1 - 2k}}\big)^{2/3}.
$$
This implies that $\#S_n(k) \tilde{=} \big({{n+1}\over{n + 1 - 2k}}\big)^{2/3}$; that is, the number of $j\in S_n(k)$ is in the range of some constant multiple of $\big({{n+1}\over{n + 1 - 2k}}\big)^{2/3}$. Since the points $\{\theta_{n,j}\}$ are evenly distributed in $S_n(k)$, we may conclude that
$$
\sum_{j\in S_n(k)} {{1}\over{n|e^{i\theta_{n,j}} - e^{i J(\theta_{n,k})}| }} \ge c_1 \log\big({{n+1}\over{n + 1 - 2k}}\big)^{2/3},
$$
for some constant $c_1 >0$. Now, in view of
$$
|1 + e^{i\theta_{n,k}}| \tilde{=}  {{n + 1 - 2k}\over{n+1}},
$$
we have $|1 - e^{i J(\theta_{n,k})}| \tilde{=} |J(\theta_{n,k})|$,\, and hence, it follows from ($\ref{ineq:theta^*00}$) that
$$
nJ(\theta_{n,k})|1 - e^{i J(\theta_{n,k})}|^{1/2}|1 + e^{i\theta_{n,k}}|^{1/2} \tilde{=} n|J(\theta_{n,k})|^{3/2}|1 + e^{i\theta_{n,k}}|^{1/2}\tilde{=} {{n + 1 - 2k}},
$$
and from ($\ref{ineq:L_nOmega2}$) that
\begin{eqnarray*}\label{ineq:L_nOmega3}
 \sum_{j = 0}^n \big | {{{\omega}_n(z_0)} \over{{\omega}'_n(z^*_{n,j})(z_0 - z^*_{n,j}) }}\big | &\ge &c_1 \sum_{k ={\lfloor{5n/12}\rfloor}}^{{\lfloor{n/2}\rfloor}}{{1}\over{n + 1 - 2k}}\log\big({{n+1}\over{n + 1 - 2k}}\big)^{2/3}\\
 &\tilde{=}& \int_{1}^{{n+1}\over12} {{\log({{n+1}\over{x}})}\over{x}} dx \\
 &= &{{\log^2(n+1) - \log^2(12)}\over{2}}.
\end{eqnarray*}
This implies that ($\ref{ineq:Lebesg0}$) in Theorem $\ref{theorem:JustFejer}$ holds.\\

In order to prove ($\ref{ineq:LebesgInfty}$), let $n=n_j$ satisfy ($\ref{ineq:closeToJ_theta}$) in Lemma $\ref{lemma:closeToJ_theta}$; that is,
$$
|J(\theta_{n_j,{\lfloor{n_j/2}\rfloor}} )- \theta_{n_j,j}| \le \frac{c_2}{n_j^{4/3}}.
$$
For $z_1 = \psi_0(e^{i(\theta_{n_j,j} + 1/n)})$ and $t_1 = \theta_{n_j,j} + 1/n$,  ($\ref{ineq:Key0}$) and ($\ref{ineq:Key'01}$) in Lemma $\ref{lemma:Omega_ToOmega1}$ imply that
$$
\big | {{{\omega}_{n_j}(z_1)}  }\big | \tilde{=} \big |{{ (z_1 - z^*_{n_j,j}) ( z - z^*_{n_j,{\lfloor{n_j/2}\rfloor}}) } \over {(z_1 - \zeta^*_{n,k_1(t_1)}) ( z_1 - \zeta^*_{n_j,{\lfloor{n_j/2}\rfloor}})}}\big|
$$
and
$$
|\omega'_n(z^*_{n_j,j})| \tilde{=} \big |{{{ z^*_{n_j,j} - z^*_{n_j,{\lfloor{n_j/2}\rfloor}} }} \over {(z^*_{n_j,j} - \zeta^*_{n_j,j}) ( z^*_{n_j,j} - \zeta^*_{n,{\lfloor{n_j/2}\rfloor}})}}\big| 
\tilde{=} \big |{{{ z^*_{n_j,j} - z^*_{n_j,{\lfloor{n_j/2}\rfloor}} }} \over {(z_1 - \zeta^*_{n_j,j}) ( z_1 - \zeta^*_{n,{\lfloor{n_j/2}\rfloor}})}}\big|.
$$
Thus, we have
\begin{eqnarray*}
\big | {{{\omega}_n(z_1)} \over{{\omega}'_n(z^*_{n_j,j})(z_1 - z^*_{n_j,j}) }}\big | &\tilde{=}& \big |{{ (z_1 - z^*_{n_j,j}) ( z_1 - z^*_{n_j,{\lfloor{n_j/2}\rfloor}}) }\over {({ z^*_{n_j,j} - z^*_{n_j,{\lfloor{n_j/2}\rfloor}})(z_1 - z^*_{n_j,j})  }}}\big |\\
&= &\big |{{   z_1 - z^*_{n_j,{\lfloor{n_j/2}\rfloor}} }\over {{ z^*_{n_j,j} - z^*_{n_j,{\lfloor{n_j/2}\rfloor}}  }}}\big |\\
&\tilde{=}& \big |{{(\theta_{n_j,j} + 1/n_j - J(\theta_{n_j,{\lfloor{n_j/2}\rfloor}}))\psi_0'(e^{i J(\theta_{n_j,{\lfloor{n_j/2}\rfloor}})})}
\over{(\theta_{n_j,j}  - J(\theta_{n_j,{\lfloor{n_j/2}\rfloor}}))\psi_0'(e^{i J(\theta_{n_j,{\lfloor{n_j/2}\rfloor}})})}} \big| \\
&=& \big |1 + {{1}
\over{n_j(\theta_{n_j,j}  - J(\theta_{n_j,{\lfloor{n_j/2}\rfloor}}))}} \big|. 
\end{eqnarray*}
By ($\ref{ineq:closeToJ_theta}$) in Lemma $\ref{lemma:closeToJ_theta}$, we obtain
$$
L_{{\bf {z^*_{n_j}}}} \ge \big | {{{\omega}_n(z_1)} \over{{\omega}'_n(z^*_{n_j,j})(z_1 - z^*_{n_j,j}) }}\big | \ge c_1\big ( 1 + \frac{{n_j^{4/3}}}{{c_2n_j}} \big ) > \frac{c_1}{c_2} \sqrt[3]{n_j}.
$$
Therefore, ($\ref{ineq:LebesgInfty}$) in Theorem $\ref{theorem:JustFejer}$ also holds.
\end{Proof}
\vskip.1in
\noindent Next we will establish {\bf Theorem 2}
\vskip .1 in
\begin{Proof}
From the distance estimation ($\ref{ineq:dist21}$), we set $\{\psi_0(e^{i\theta_{n,k}}): |\theta_{n,k}| \le \frac{2\pi}{3}\}$ which satisfies the separating condition ($\ref{ineq:sep_c0}$):
$$
\min_{ k \neq j, |\theta_{n,k}| \le \frac{2\pi}{3}, |\theta_{n,j}| \le \frac{2\pi}{3} }   {{|z^*_{n,j}-z^*_{n,k}|} \over {\min \big(\text{dist}(z^*_{n,j},\Gamma_n ),\text{dist}(z^*_{n,k},\Gamma_n ) \big)}} \ge c_1>0.
$$
By Lemma 2.1 and Lemma 2.2 in \cite{CZ2021}, we have
$$
c_1 \sum_{|\theta_{n,k}| \le \frac{2\pi}{3}} |P_n (z^*_{n,k} ) |^p \text{dist}(z^*_{n,k},\Gamma_n) \le \int_\gamma |P_n (z) |^p |dz|.
$$
Similarly, we also have
$$
c_1 \sum_{|\theta_{n,k}| > \frac{2\pi}{3}} |P_n (z^*_{n,k} ) |^p \text{dist}(z^*_{n,k},\Gamma_n)  \le \int_\gamma |P_n (z) |^p |dz|.
$$
Hence, for all $P_n(z)\in \Pi_n$, 
\begin{equation} \label{ineq:MzLeft}
c_1 \sum^n_{k = 0} |P_n (z^*_{n,k} ) |^p \text{dist}(z^*_{n,k},\Gamma_n) \le 2\int_\gamma |P_n (z) |^p |dz|.
\end{equation}
Suppose, on the contrary, that Theorem 2 is false; that is,
\begin{equation} \label{ineq:MzRight}
\lim\sup_{n\rightarrow\infty} \max_{P_n(z)\in\Pi_n} \frac{\int_{\gamma_0} |P_n(z)|^p |dz| }{\sum^n_{k = 0} |P_n(z^*_{n,k})|^p \text{dist}(z^*_{n,k},\Gamma_n)} < \infty,
\end{equation}
 for some sub-sequence $\{n_{\ell}\}\rightarrow\infty$. Then it follows from ($\ref{ineq:MzLeft}$) that the modified Marcinkiewicz-Zygmund inequalities are valid for $z_{n_{\ell},k} = z^*_{n_{\ell},k}$ in ($\ref{mod MZ ineq}$). Therefore, by a result in \cite{CZ2021} (see Theorem 2.1 of \cite{CZ2021}), the separate condition ($\ref{ineq:sep_c0}$) for the point $z^*_{n, {\lfloor{n/2}\rfloor}}$, with $n=n_{\ell}$ is satisfied. However, from the distance estimation ($\ref{ineq:dist_00}$) for $\theta_{n,{\lfloor{n/2}\rfloor}}$, we have
$$
\min_{|\theta_{n,j}|<{{2\pi}\over{3}}} (n+1)^{\frac{4}{3}} |\theta_{n,j} - J(\theta_{n,{\lfloor{n/2}\rfloor}})| \ge  {{c_1}\sqrt[3]{n+1}}, 
$$ 
again for $n=n_{\ell}$. This is a contradiction to ($\ref{ineq:closeToJ_theta}$) of Lemma $\ref{lemma:closeToJ_theta}$; hence, completing the proof of the theorem.
\end{Proof}
\vskip.1in

\noindent Finally we now establish {\bf Theorem 3}.
\vskip.1in
\begin{Proof}
The proof of the first inequality in  ($\ref{ineq:Lebesg2}$) is similar to that of ($\ref{ineq:Lebesg0}$) in Theorem $\ref{theorem:JustFejer}$. For the proof of the second inequality, let $z = \psi_0(e^{it})$, with $0\le t\le {{2\pi}\over{3}}$. Then from the first item  of ($\ref{equivalence}$) in Lemma $\ref{lemma:Key00}$, we have
$$
|\tilde\omega(z)| \le c_2;
$$
and under the condition of $0\le\theta_{n,j}\le {{2\pi}\over{3}}$, it follows from ($\ref{equivalence}$) in Lemma $\ref{lemma:Key00}$, with the simplified notation $k = k_2(\tilde\theta^*_{n,j})$ that
\begin{equation}\label{ineq:Key'14}
|\tilde\omega'_n(z^{**}_{n,j})| \tilde{=}
\big |{{ z^{**}_{n,j} - z^{**}_{n,k} } \over {(z^*_{n,j} - \zeta^*_{n,j}) ( z^*_{n,k} - \zeta^*_{n,k})}}\big|.
\end{equation}
Hence, by Lemma $\ref{lemma:localSum}$, we may conclude that 
\begin{eqnarray} \label{sumSeparateCloseToT}
\sum_{\{j:0\le\theta_{n,j}\le {{2\pi}\over{3}}\}} \big | {{{\tilde\omega}_n(z)} \over{\tilde\omega'_n(z^{**}_{n,j})(z - z^{**}_{n,j}) }}\big | &\le&  \sum_{\{j:|k_2(\theta_{n,j}) - k_2(t)| \le 1\}} \big | {{\tilde\omega(z)} \over{\tilde\omega'_n(z^{**}_{n,j})(z - z^{**}_{n,j}) }}\big | \nonumber\\
&+& \sum_{\{j:|k_2(\theta_{n,j}) - k_2(t)| > 1\}} \big | {{c_2} \over{\tilde\omega'_n(z^{**}_{n,j})(z - z^{**}_{n,j}) }}\big | \nonumber\\
& \le &  c_2 \log(n) + c_2\sum_{\{j:|k_2(\theta_{n,j}) - k_2(t)| > 1\}} \big | {{(z^*_{n,j} - \zeta^*_{n,j}) ( z^*_{n,k} - \zeta^*_{n,k})} \over{(z^{**}_{n,j} - z^{**}_{n,k}) (z - z^{**}_{n,j}) }}\big |.
\end{eqnarray}
Also, for ${{5\pi}\over{6}}\le\theta_{n,k} < \pi$, we have
$$
|z - z^{**}_{n,j}| = |\psi_0(e^{i J^{-1}(t)}) - \psi_0(e^{iJ^{-1}(\tilde\theta^*_{n,j})})| \tilde{=} |J^{-1}(t) - J^{-1}(\tilde\theta^*_{n,j})|(\pi - J^{-1}(t) + |J^{-1}(t) - J^{-1}(\tilde\theta^*_{n,j})|)^{-1/2},
$$
and under the condition of $|k - k_2(t)| = |k_2(\tilde\theta^*_{n,j}) - k_2(t)| > 1$, it follows that
$$
|J^{-1}(t) - J^{-1}(\tilde\theta^*_{n,j})| \tilde{=} |J^{-1}(t) - \theta_{n,k}|.
$$
Therefore,
$$
|z - z^{**}_{n,j}| \tilde{=} |J^{-1}(t) - \theta_{n,k}|(\pi - J^{-1}(t) + |J^{-1}(t) - \theta_{n,k}|)^{-1/2} \tilde{=} |J^{-1}(t) - \theta_{n,k}|(\pi - \theta_{n,k} + |J^{-1}(t) - \theta_{n,k}|)^{-1/2},
$$
so that
$$
{{|z^*_{n,k} - \zeta^*_{n,k}|}\over{|z - z^{**}_{n,j}|}} \tilde{=} {{|\pi - \theta_{n,k}|^{-1/2}}\over{n|z - z^{**}_{n,j}|}} \tilde{=} {{(\pi - \theta_{n,k} + |J^{-1}(t) - \theta_{n,k}|)^{1/2}}\over{n| J^{-1}(t) - \theta_{n,k}||\pi - \theta_{n,k}|^{1/2}}},
$$
and
\begin{equation}\label{est11}
{{|z^*_{n,k} - \zeta^*_{n,k}|}\over{|z - z^{**}_{n,j}|}} \le
{{c_2}\over{n| J^{-1}(t) - \theta_{n,k}|}} + {{c_2}\over{n| J^{-1}(t) - \theta_{n,k}|^{1/2}|\pi - \theta_{n,k}|^{1/2}}}, \hskip15pt |k_2(\tilde\theta^*_{n,j}) - k_2(t)| > 1.
\end{equation}
Next, for ${{2\pi}\over{3}}<\theta_{n,k} <{{5\pi}\over{6}}$, we have
$$
|z - z^{**}_{n,j}|\tilde{=} | J^{-1}(t) - J^{-1}(\tilde\theta^*_{n,j})|(|{{2\pi}\over{3}} - J^{-1}(\tilde\theta^*_{n,j})| + | J^{-1}(t) - J^{-1}(\tilde\theta^*_{n,j})|) \ge c_1 | J^{-1}(t) - \theta_{n,k}|(|{{2\pi}\over{3}} - \theta_{n,k}| + {{2\pi}\over{n+1}}).
$$
This yields
$$
{{|z^*_{n,k} - \zeta^*_{n,k}|}\over{|z - z^{**}_{n,j}|}} \tilde{=} {{(|{{2\pi}\over{3}} - \theta_{n,k}| + 1/n)}\over{n|z - z^{**}_{n,j}|}} \le {{c_2}\over{n| J^{-1}(t) - \theta_{n,k}|}},
$$
so that ($\ref{est11}$) holds for ${{2\pi}\over{3}}\le\theta_{n,k} < \pi$. In the consideration of the summation in left-hand side of (\ref{sumSeparateCloseToT}), we observe that

\begin{eqnarray*}
& &\sum_{\{j:|k_2(\theta_{n,j}) - k_2(t)| > 1\}} \big | {{(z^*_{n,j} - \zeta^*_{n,j}) ( z^*_{n,k} - \zeta^*_{n,k})} \over{(z^{**}_{n,j} - z^{**}_{n,k}) (z - z^{**}_{n,j}) }}\big |\\
&\le & \sum_{{{2\pi}\over{3}}\le\theta_{n,k} < \pi,|\theta_{n,k} -  J^{-1}(t)|>1/n}\big({{c_2}\over{n| J^{-1}(t) - \theta_{n,k}|}} + {{c_2}\over{n| J^{-1}(t) - \theta_{n,k}|^{1/2}|\pi - \theta_{n,k}|^{1/2}}}\big)  \sum_{j\in S_n(k)}\big | {{z^*_{n,j} - \zeta^*_{n,j} } \over{z^{**}_{n,j} - z^{**}_{n,k} }}\big |\\
&\le & \sum_{{{2\pi}\over{3}}\le\theta_{n,k} < \pi,|\theta_{n,k} -  J^{-1}(t)|>1/n}\big({{c_2}\over{n| J^{-1}(t) - \theta_{n,k}|}} + {{c_2}\over{n| J^{-1}(t) - \theta_{n,k}|^{1/2}|\pi - \theta_{n,k}|^{1/2}}}\big)\,\, \log(n),
\end{eqnarray*}
which implies that
\begin{equation}
\sum_{0\le\theta_{n,j}\le {{2\pi}\over{3}}} \big | {{{\tilde\omega}_n(z)} \over{\tilde\omega'_n(z^{**}_{n,j})(z - z^{**}_{n,j}) }}\big | \le c_2\log^2(n).
\end{equation}

As to the summation over ${{2\pi}\over{3}}<\theta_{n,j}<\pi$, it follows from ($\ref{equivalence}$) in Lemma $\ref{lemma:Key00}$, that
\begin{eqnarray*}
\big|{{\tilde\omega_n(z)}\over{{\tilde\omega'_n(z^{**}_{n,j}) (z - z^{**}_{n,j}) }   }}\big| &\tilde{=}&
\big|{{ (z - z^{**}_{n,k_1(t)}) ( z - z^{**}_{n,k_2(t)})(z^{**}_{n,j} - \zeta^*_{n,j})  } \over {(z - \zeta^*_{n,k_1(t)}) ( z - \zeta^*_{n,k_2(t)})(z - z^{**}_{n,j})}}\big|\\
&\le&
c_2\,\big|{{  ( z - z^{**}_{n,k_2(t)})(z^{**}_{n,j} - \zeta^*_{n,j})  } \over { ( z - \zeta^*_{n,k_2(t)})(z - z^{**}_{n,j})}}\big|.\\
\end{eqnarray*}
Therefore, 
$$
\sum_{{{2\pi}\over{3}}<\theta_{n,j}<\pi} \big | {{{\tilde\omega}_n(z)} \over{\tilde\omega'_n(z^{**}_{n,j})(z - z^{**}_{n,j}) }}\big |
\le c_2 \sum_{{{2\pi}\over{3}}<\theta_{n,j}<\pi, j\neq k_2(t)}\big|{{ z^{**}_{n,j} - \zeta^*_{n,j}  } \over { z - z^{**}_{n,j}}}\big| + c_2 |{{ z^{**}_{n,k_2(t)} - \zeta^*_{n,k_2(t)}  } \over { z - \zeta^*_{n,k_2(t)}}}\big|.
$$
Similar to ($\ref{est11}$), we may obtain
\begin{equation}\label{est12}
{{|z^*_{n,j} - \zeta^*_{n,j}|}\over{|z - z^{**}_{n,j}|}} \le
{{c_2}\over{n| J^{-1}(t) - \theta_{n,j}|}} + {{c_2}\over{n| J^{-1}(t) - \theta_{n,j}|^{1/2}|\pi - \theta_{n,j}|^{1/2}}}, \hskip15pt j\neq k_2(t), |\theta_{n,j}| > {{2\pi}\over{3}}.
\end{equation}
Therefore, we have
\begin{eqnarray*}
& &\sum_{{{2\pi}\over{3}}<\theta_{n,j}<\pi} \big | {{{\tilde\omega}_n(z)} \over{\tilde\omega'_n(z^{**}_{n,j})(z - z^{**}_{n,j}) }}\big | \\
&\le &  \sum_{{{2\pi}\over{3}}<\theta_{n,j}<\pi, j\neq k_2(t)}\big({{c_2}\over{n| J^{-1}(t) - \theta_{n,j}|}} + {{c_2}\over{n| J^{-1}(t) - \theta_{n,j}|^{1/2}|\pi - \theta_{n,j}|^{1/2}}}\big) + c_2 \\
&\le& c_2 \,\,\log(n).
\end{eqnarray*}
which yields
\begin{equation}\label{sum000}
\sum_{0\le\theta_{n,j}<\pi} \big | {{{\tilde\omega}_n(z)} \over{\tilde\omega'_n(z^{**}_{n,j})(z - z^{**}_{n,j}) }}\big | \le c_2\,\log^2(n).
\end{equation}
In the summation for $\theta_{n,j}<0$, let $\tilde{n} = 2\lfloor{n/2}\rfloor + 1$. Then
$$
\big|{{{\tilde\omega}_n(z)} \over{\tilde\omega'_n(z^{**}_{n,j})(z - z^{**}_{n,j}) }}\big | =
{{{|\tilde\omega}_n(z)|} \over{|\tilde\omega'_n(z^{**}_{n,\tilde{n} - j})||z - \overline z^{**}_{n,\tilde{n} -j}| }}\le {{{|\tilde\omega}_n(z)|} \over{|\tilde\omega'_n(z^{**}_{n,\tilde{n} -j})||z - z^{**}_{n,\tilde{n} -j}| }}\,\,,
$$
so that
\begin{equation}
\sum_{-\pi<\theta_{n,j}<0} \big | {{{\tilde\omega}_n(z)} \over{\tilde\omega'_n(z^{**}_{n,j})(z - z^{**}_{n,j}) }}\big | \le\sum_{0<\theta_{n,j}<\pi} \big | {{{\tilde\omega}_n(z)} \over{\tilde\omega'_n(z^{**}_{n,j})(z - z^{**}_{n,j}) }}\big |\le c_2\log^2(n).
\end{equation}
Combining this with $(\ref{sum000})$, we may conclude that second inequality of $(\ref{ineq:Lebesg2})$ in Theorem $\ref{theorem:AdJustFejer}$ also holds. This completes the proof of the theorem.
\end{Proof}

\bhag{Problems and results inspired by the work of Prof. Korevaar in the decade of 1960's}

 Although the study of Lagrange polynomial interpolation at the Fej\'er points ${\bf z}_n^{*}$ and adjusted Fej\'er points ${\bf z}_n^{**}$ on an open arc $\gamma\subset \CC$ is relatively complete in this paper, yet it has raised more questions than it answers. This section is devoted to the discussion of the notion of Lebesgue constants for more general polynomial approximation in the complex plane, the limitation of exterior conformal mapping for the selection of polynomial interpolation nodes on open arcs; and most importantly, a somewhat extensive discussion of the problems and results inspired by the work of Prof. Korevaar in the decade of 1960's, on approximation by polynomials in the complex plane with restricted zeros and asymptotically neutral distribution of electrons, with emphasis on the placement of point-masses in the Euclidean space $\RR^s$, for $s\ge 2$, and the inverse problem of locating such point-masses in terms of the ``super-resolution" recovery problem.

\subsection{Lebesgue constants and limitation of exterior conformal mapping for interpolation node selection on open arcs in the complex plane} 

Let $\Pi_n$ denote the space of polynomials of degree $\le n$; and for a bounded set $X\subset \CC$, let $B_X$ be a Banach space of functions defined on $X$, such that $\Pi_n \subset B_X$. In this paper, for a bounded linear operator $A_{n}$ from $B_X$ to $\Pi_n$ that possesses the polynomial preservation property, namely: $A_{n}p_n\, = \,p_n$ for all $p_n \in \Pi_n$, the ``Lebesgue constant" corresponding to $A_{n}$ is defined by the operator norm $\|A_{n}\|$. In Lemma 1.1, we have verified that the operator norm of the linear operator $A_{{\bf {z}}_n}$ defined in ($\ref{operator}$) agrees with the commonly used definition of Lebesgue constant $L_{{\bf {z}}_n}= \max_{z\in X} \, \sum^{n}_{j=0} |\ell_ {n,j}(z)|$ as shown in (\ref{Lcformulation}). However, there are other definitions of Lebesgue constants in the literature that require verification to agree with the above operator-norm definition. In particular, for a Jordan domain $D\subset \CC$ with boundary being a simple closed piece-wise smooth curve $\Gamma$, we may consider the (generalized) Hardy space $H^{p}(D)$, for $p\ge 1$, of analytic functions $f$ in $D$ with non-tangential limit $f^{*} \in L_p(\Gamma)$, as our Banach space $B(X)$ of functions $f^{*}$, defined (almost everywhere) on $X=\Gamma$, and the linear operator $A_{F_n}$, that maps $H^{p}(D)$ to the $n$-th partial sum $S_n(\cdot; f)$ of the Faber series representation of $f$ in $D$. Then T. K\"ovari and C. Pommerenke \cite{KOV1966} introduced the following Lebesgue constant:
\begin{equation} \label{KP_Lebesgue}
L_{n,D} := \max_{z\in\Gamma}\int_{0}^{1}\big|\sum_{k = 0}^n F_{k}(z) \exp\{-i(2(k+1)\pi\omega\}\big | d\omega,
\end{equation}
where $\{F_k(z)\}_{k=0}^n$ are the Faber polynomials for the domain $D$. It is clear that $A_{F_n}p_n\, = \,p_n$ for all $p_n \in \Pi_n$, but it is not clear if $L_{n,D} = A_{F_n}$. 

In this paper, we consider polynomial interpolation at the Fej\'er points on open arcs $\gamma$ in the complex plane, and discover the limitation of the exterior conformal map from $|w|>1$ to $\CC^{*}\backslash \gamma$ for this problem. While the conformal map is a very powerful tool in complex analysis, it happens that its  ``angle-preserving" property for all $z\in \CC\backslash \gamma$, when extending to the boundary, imposes limitation of its application to interpolation node selection for open arcs $\gamma$ with corners, for the reason that not only the tangent at $z\in\gamma$ is not unique at the corners of $\gamma$, but also the images of the (rotated) roots of unity of the conformal mapping of the level curves $|w|= 1+\frac{1}{m}$ to ``both sides" of $\gamma$ that converge to the Fej\'er points on $\gamma$ could be very very close to one another, particularly near the corners, when $m\rightarrow {\infty}$. For this reason, a procedure of adjusting the Fej\'er points proposed in Section 2 is necessary to obtain more precise rate of growth of the Lebesgue constant for Lagrange interpolation at Fej\'er points. However, we believe that our adjustment procedure is not natural, leading to our discussion of the earlier work of J. Korevaar on the distribution of polynomial approximation with restricted zeros in the next sub-section, and the placement of point-masses in $\RR^s$ for $s\ge 1$, such as ``electrons", in the final sub-section. 

\subsection{Polynomial approximation with restricted zeros, distribution of electrons, and total energy of electrostatic fields}

 During the decade of the 1960's, inspired by the paper \cite{GM} of G.R. MacLane, Professor J. Korevaar has invested a vast amount of time and energy on the research topics of ``polynomial approximation on a Jordan domain $D \subset \CC$ with restricted zeros" and the related problem of ``asymtotically neutral distributions of electrons". His contributions in these directions include his pioneering paper \cite{K1}, the other work \cite{K2, K3} and his joint papers \cite{KC, KM1, KM2}, as well as four of his students: M.D. Thompson, J.M. Elkins, C.K. Chui, and D.T. Piele (see portions of their Ph.D. theses published in \cite{T1}, \cite{JE}, \cite{C,C1}, and \cite{P1, P2}, respectively). Since  a polynomial $P_n\in \Pi_{n+1}$ with zeros ${\bf z}_n =\{z_{n,k}\}_{k=0}^n\subset\Gamma$, where $\Gamma$ denotes the boundary of the Jordan domain $D$, can be written as $P_{n+1}(z) = c_0\prod_{k = 0}^n(z - z_{n,k})$ for some $c_0\neq 0$, its logarithmic derivative becomes
\begin{equation} \label{elec field}
F_{{{\bf z}_n}} (z):=\frac{d}{dz}\log P_{n+1}(z)= \sum_{k = 0}^n \frac{1}{z - z_{n,k}},
\end{equation}
which represents the electrostatic field at $z\in D$ due to the electrons located at the positions $z_{n,k}\in \Gamma$.
In \cite{K1}, a family ${\bf z}_n =\{z_{n,k}\}_{k=0}^n$ of points on $\Gamma$ is defined to be ``asymptotically neutral", if $F_{{{\bf z}_n}} (z)\rightarrow 0$ uniformly on every compact subset of $D$, for $n=n_j \rightarrow \infty$. It was also shown in 
(Lemma 4.1 of \cite{K1}) that the family $\{z_{n,k} = \psi_{\Gamma}(e^{i(\theta +\frac{2\pi k}{n+1})})\}$ of  Fej\'er points is asymptotically neutral for almost all $\theta$, where $\psi_\Gamma$ denotes the exterior conformal map from $|w|>1$ to $\CC^{*}\backslash D$, as considered in this present paper.

This concept of ``asymptotically neutral family of points" introduced in Korevaar's paper \cite{K1} has significantly high impact to the past and current research areas. Firstly, while the family $\{z_{n,k} = \psi_{\Gamma}(e^{i(\theta +\frac{2\pi k}{n+1})})\}$ is asymptotically neutral for almost all $\theta$, the nature of the corresponding ``total energy" of  the electrostatic field in $D$, defined by 
\begin{equation} \label{energy}
E({\bf z}_n; D) := \int_{D}|F_{{\bf z}_n} (x+iy)|dx\,dy,
\end{equation}
is still not completely well understood, where $D$ is considered as a domain in $\RR^2$. In this regard, to respond to the conjecture posed by the first author in \cite{C5} that for the unit disc $D_0 := \{z: |z|<1\}$, the distribution of electrons that gives rise to minimum total energy $E({\bf z}_n; D_0)$, defined in (\ref{energy}), is attained by ${\bf z}_n= {\bf z}^{*}_n$ for
$${\bf z}^{*}_n := \{e^{i(\theta +\frac{2\pi k}{n+1})}\}_{k=0}^{n}\,,$$ 
for any $\theta\in\RR$, D.J. Newman proved in \cite{DN} that  $E({\bf z}_n; D_0)$ is bounded below by $\frac{\pi}{18}$ for all choices of ${\bf z}_n$ on the unit circle $|z|=1$ and for all $n\in \ZZ_{+}$. In a more recent paper \cite{ABF}, Abakumov, Borichev and Fedorovskiy modified this conjecture in the context of weighted Bergman spaces of square-integrable functions (in Theorem 1 of \cite{ABF}), and proved that the minimum (Bergman-modified) energy $F^{B}_{{\bf z}_n; D_0}(z)$, for ${\bf z}_n$ on the unit circle is indeed attained at (arbitrary rotation of) the $(n+1)$-st roots of unity for all $n\in \ZZ_{+}$. However, to the best of our knowledge, the original conjecture posed in \cite{C5} is still open. In this regard, we wonder if Newman's result in \cite{DN} remains valid for more general Jordan domains $D\subset\CC$ with simple closed boundary curves $\Gamma$. Of course, the conjecture that the minimum total energy $E({\bf z}_n; D)$ is attained by the Fej\'er points on $\Gamma$ is a much harder problem.

On the other hand, concerning the approximation of analytic functions in the Jordan domain $D\subset\CC$ by functions $F_{{{\bf z}_n}} (z)$ defined in (\ref{elec field}), for $\{{\bf z}_{n}\} \subset \Gamma$ studied by Korevaar in \cite{K1}, the results for $D$ being the special case of the unit disc $D_0$ include the paper \cite{T2} of M. Thompson on uniform approximation on compact subsets of the unit disk $D_0$ of functions of $H^{\infty}$ with poles ${\bf z}_{n}$ on the unit circle, and the recent paper \cite{ABF} of Abakumov, Borichev and Fedorovskiy quoted above, with details as well as precise orders of approximation and asymptotic behaviour, in the (weighted) Hilbert-Bergman spaces, for a wide class of weights. As to the general Jordan domain $D\subset\CC$, approximation the Bers space by functions $F_{{{\bf z}_n}} (z)$, with poles ${\bf z}_n$ restricted to the simple closed boundary curve $\Gamma$ of $D$, was first studied in the paper \cite{C6}, and the precise orders of approximation was later established in the joint paper \cite{CS} of X. Shen and the first author of the present paper. 

The results on approximation of analytic function in a simply connected bounded domain $D\subset\CC$ by polynomials with zeros that lie on the boundary curve $\Gamma$ of $D$ naturally carry over to the possibility of approximation of (real-valued) harmonic functions in $D$ by functions $u_n(z)$ of the form:

\begin{equation} \label{harmonic 1}
u_n(z) = c_n + \sum_{k=0}^{n} \log|z-z_{n,k}|,\,\,\,\,c_n\in\RR,\,\,\,z_{n,k}\in \Gamma\,,
\end{equation}
for $z\in D\subset\CC$, in that every real-valued function $U(z)$, harmonic in $D$, can be approximated, uniformly on compact subsets of D, by $u_n(z)$ in (\ref{harmonic 1}) for some sequence of real numbers $c_n$. We say that the harmonic function $U(z)$ is $\Gamma$-approximable. The reason for introducing the notion of ``$\Gamma$-approximable" harmonic functions in the joint paper \cite{KC} by Korevaar and the first author of this paper is that for $\Gamma=\Gamma_1 \cup\cdots\cup \Gamma_m$, where $\Gamma_1,\cdots, \Gamma_m$ are analytic boundary curves of $m$ mutually disjoint bounded Jordan domains $D_1,\cdots,D_m$, respectively for $m\ge2$, then the class of ``$\Gamma$-approximable" harmonic functions depends on some conditions of the harmonic measures of the boundary curves at $z=\infty$. While the results in \cite{KC} are valid for arbitrary $m\ge2$, for convenience we only consider $m=2$ in this discussion. Let $E:= \CC^{*}\backslash \text{clos}\,{(D_1\cup D_2)}$ and $\omega_1(z)$ denote the harmonic measure of $\Gamma_1$ relative to $E$. Then the harmonic measure of $\Gamma_2$ relative to $E$ is $\omega_2(z)= 1-\omega_1(z)$. Also, let $\Omega_1 = \omega_1(\infty)$. In \cite{KC}, it is shown that if  $\Omega_1$ is an irrational number, then all harmonic functions $U(z)$ in $D=D_1\cup D_2$ are $\Gamma$-approximable. On the other hand, if $$\Omega_1= \frac{p}{q},$$ where $p$ and $q$ are relatively prime positive integers, then a function $U(z)$ harmonic in $\text{clos}\, (D)$ is $\Gamma$-approximable, if and only if it satisfies:
\begin{equation} \label{approximable}
\frac{1}{2\pi}\int_{\partial D} U(z)\frac{\partial \omega_1(z)}{\partial N} |dz| \equiv 0 \big(\text{mod}\,\frac{1}{q}\big),
\end{equation}
where $N$ denotes the outer normal to $\partial D$.

The study of approximation of harmonic functions by $u_n(z) = c_n + \sum_{k=0}^{n} \log|z-z_{n,k}|$ in $D$, introduced in (\ref{harmonic 1}), with $z_{n,k}\in \partial D$, was extended to the $3$-dimensional Euclidean space $\RR^3$ by Korevaar's former student, D.T. Piele in \cite{P2}, where $\log|z-z_{n,k}|$ is replaced naturally by $|{\bf x} - {\bf x}_{n,k}|^{-1}$ for bounded open connected domain $D$ with connected complement $\RR^3\backslash D$, where $|{\bf x}|$ denotes the Euclidean norm of ${\bf x}\in \RR^3$. More precisely, it is proved in \cite{P2} that, analogous to (\ref{harmonic 1}), every harmonic function $U({\bf \x})$ in $D\subset \RR^3$ can be approximated, uniformly on compact subsets of $D$ by 
\begin{equation} \label{harmonic 2}
u_n({\bf \x}) = C_n + \sum_{k=0}^{n} \frac{1}{|{\bf x} - {\bf x}_{n,k}|},\,\,\,\,C_n\in\RR,\,\,\,{\bf x}_{n,k}\in \partial D\,,
\end{equation}
for ${\bf x}\in D$ and some subsequence $\{u_{n_j}({\bf \x})\}$ of $\{u_n({\bf \x})\}$.

It is interesting to investigate the extension from $D\subset \RR^3$ to mutually disjoint unions of bounded open connected domains $D_1,\cdots,D_m$ in $\RR^3$, and if the notion of $\partial D$-approximable harmonic functions and an analogy of the condition on the ``harmonic measure at $\infty$" as in \cite{KC} for the extension make sense.

\subsection{Placement of point-masses in higher dimensions and super-resolution point-mass recovery}

Another problem area of great interest to Prof. Korevaar is the placement of electrons on the boundary $\partial D$ of a bounded open connected domain $D\subset\RR^s$ for $s\ge 2$ to minimize the electrostatic field in $D$. Of course, this problem is more meaningful for the 3-dimensional Euclidean space $\RR^3$, and in the joint paper \cite{KM1} with J.L.H. Meyers, he considered the ``spherical Faraday cage" and introduced a Chebyshev-type quadrature on the unit sphere $S^2\subset \RR^3$, as the boundary of the unit ball $D_0$. Since the formula of the electrostatic field $F_{{\bf z}_n} (z)$ for $z\in \CC$ due to the electrons ${\bf z}_n$ in (\ref{elec field}) can be written as the complex conjugate of $F_{{{\bf z}_n}} (z)$, namely:
$$ \sum_{k = 0}^n \frac{z - z_{n,k}}{|z - z_{n,k}|^2}\,\,,$$
it is natural to define the ``electrostatic field" in $\RR^s$ for $s\ge3$ by
\begin{equation} \label{elec fieldR3}
\F_{{{\bf X}_n}} ({\bf x})= \sum_{k = 0}^n \frac{\bf x - {\bf x}_{n,k}}{|{\bf x} - {\bf x}_{n,k}|^s}\,\,,
\end{equation}
where ${\bf x}$ and ${{\bf X}_n}=\{{\bf x}_{n,k}\}_{k=0}^n$ are in $\RR^s$.  

In considering the ``spherical Faraday cage" $S^2:=\{{\bf x}:|{\bf x}| = 1\} \subset \RR^3$, Prof. Korevaar was interested in the problem of placing the electrons (called point-mass with all coefficients equal to 1, to be discussed later) at the Fekete points ${\bf X}^{*}_n =\{{\bf x}^{*}_{n,k}\}_{k=1}^n$ on the unit sphere $S^2$, defined by 
\begin{equation} \label{fekete}
\{{\bf x}^{*}_{n,k}\}:= \arg \,\max\, \prod^n_{j = 0}\prod_{k = 0, k \neq j}^n |{\bf x}_{n,k} - {\bf x}_{n,j}|\,,
\end{equation}
where $|{\bf x}_{n,k}| = 1$ for $k=0,\dots,n$ and all $n \in \ZZ_{+}$. However, in the present paper, we consider the Fekete points ${\bf X}_n^{*} := \{{\bf x}^{*}_{n,k}\}_{k=0}^n$, defined in (\ref{fekete}), to be placed on the surface $\partial D$ of any bounded open connected domain $D\subset\RR^s$ for $s\ge 2$, by replacing the restriction $|{\bf x}_{n,k}| = 1$ to ${\bf x}_{n,k}\in \partial D$ in the above product.   
In addition, we extend the definition of total energy in (\ref{energy}) from the complex plane $\CC$ to $\RR^s$ for $s\ge 3$ as follows: 
\begin{equation} \label{energy2}
\E({\bf X}_n; D) := \int_{D}|\F_{{\bf X}_n} ({\bf x})|\,d{\bf x}\,\,,
\end{equation}
and propose the challenging open problem of determining the bounded open connected domains $D\subset \RR^s$ and their corresponding surfaces $\partial D$, for which the {\bf``minimum"} of the total energy $\E({\bf X}_n; D)$ is attained at the Fekete points ${\bf X}_n^{*} := \{{\bf x}^{*}_{n,k}\}_{k=0}^n \subset \partial D$. 

In this regard, we refer the interested reader to the paper \cite{LB3} by Len Bos in this special issue that provides examples with computer codes to investigate Fekete points on a  simplex in $\RR^s$ for $s\ge3$. In addition, the interested reader is also referred to the design procedures proposed by E. Saff and Kuijlaars in \cite{SK} and their paper \cite{KS} for placing the electrons on the unit sphere $S^2$. While our proposed problem is very difficult to solve, it would be fruitful to consider computer simulation to minimize the total energy for such typical domains $D$ as the standard simplex and the unit ball in $\RR^2$ and $\RR^3$ over the points ${\bf X}_n := \{{\bf x}_{n,k}\}_{k=0}^n$ on the boundary curve and boundary surface, respectively, to investigate if the minimum is attained at the Fekete points.

For real-world applications, the points ${\bf X}_n := \{{\bf x}_{n,k}\}_{k=0}^n\subset \RR^s$, where $s\ge1$, are associated with some constants $\{c_{n,k}\}_{k=0}^n\subset \RR$; and for each $n\in \ZZ_{+}$, the set of pairs $\{(c_{n,k}, {\bf x}_{n,k})\}_{k=0}^n$ is called a point-mass, represented by the distribution
\begin{equation}
\label{pointmass in s-D}
h({\bf y}) = \sum_{k=0}^{n} c_{n,k}\delta ({\bf y} - {\bf x}_{n,k}), 
\end{equation}
where $\delta$ denotes, as usual, the Dirac delta. The super-resolution problem is the inverse problem of recovering the point-mass $\{(c_{n,k}, {\bf x}_{n,k})\}_{k=0}^n$ from the Fourier transform $D({\bf x})$ of $h({\bf y})$, called the data function, namely:
\begin{equation}
\label{data fn} 
D({\bf x}) =  \sum_{k=0}^{n} c_{n,k} \exp\{-i{\bf x}_{n,k}\cdot {\bf x}\}.
\end{equation}
This problem was first addressed in the pioneering paper \cite{D} of D. Donoho for the special case of $\{{\bf x}_{n,k}\}_{k=0}^n$ that lies on some unknown rectangular grid in $\RR^2$. The general setting is formulated in \cite{C3}, where clusters of points are considered in \cite{C4} to represent the geometric shapes of objects, such as imagery of the celestial bodies captured by the James Webb space telescope, and single-strand helix images of DNA molecules and double-strand helix images of RNA molecules in a human cell, observed by the ``single molecule microscope". The current active research direction of super-resolution is considered as a research area in ``Big data". The interested reader is referred to the references in \cite{C3, C4} and there-in for the vast amount of literature in this area. In \cite{C3}, three families of super-resolution (SR) wavelets are introduced, including the Gaussian SR (GSR) wavelet $\Psi^g_{v,n}$ and the corresponding continuous wavelet transform (CWT), also called integral wavelet transform (see \cite{C2}), defined by
\begin{equation}\label{GSRW}
\begin{cases}
{\displaystyle \Psi^g_{v,n} (x) := \Big(\psi_n * g_{\kappa,v,n}\Big)(x)};\\
&\\
 {\displaystyle (W_{\Psi^g_{v,n}}D)(t, a) :=\int_{\RR}\overline{\Psi^g_{v,n} \Big(\frac{x-t}{a}\Big)}D(x)\,\frac{dx}{a}}\,\,,
\end{cases}
\end{equation}
where the wavelet $\psi_n(x)$, called the $n$-th order Haar wavelet in \cite{C3}, has vanishing moments of order $n\ge1$ and is defined by
\begin{equation}\label{Haar wavelets}
\begin{cases}
{\displaystyle \psi_n(x):= N^{(n)}_{2n}(x + n)};\\
&\\ 
    {\displaystyle N_{2n}(x):= ({\underbrace{\chi_{[0,1)} * \cdots * \chi_{[0,1)}}_{2n}})(x)},
\end{cases}
\end{equation}
and where $g_{\kappa,v,n}(x) := (-i)^n g_v(x)e^{-i\kappa x}$,\,\,with $g_v(x)$ being the Gaussian probability density function; $D(x)$ is the data function in (\ref{data fn}), and \,$\kappa \doteq 2.331122371$ is the unique solution of the equation $x= \tan\,x$, for $x>0$. By setting $t=0$ in (\ref{GSRW}), we obtain the search function $F^{g}_{v,n}(a):= (W_{\Psi^g_{v,n}}D)(0, a)$ for extracting the points $x_{n,k}$, by matching the variable $a$ of search function $F^{g}_{v,n}(a)$ with $\frac{\kappa}{x_{n,k}}$ for each $k=0, \dots,n$. In resolving an inverse problem of super-resolution, the parameter $v>0$ plays the role of narrowing the search window for isolating the points $x_{n,k}$ for increasing values of the variance of the Gaussian function $g_v(x)$, and the parameter $n$ plays the role of amplifying the coefficients $c_{n,k}$ by the multiplication factor of $\xi^{n}$, where 
\begin{equation*}
 \xi:=|\hat \psi_{1}(\kappa)|\doteq 1.449222080.    
\end{equation*}
In the above discussion, we only show the univariate GSR wavelet $\Psi^g_{v,n}$ with variance $v^2$ and only the corresponding univariate CWT. The interested reader is referred to the paper \cite{C3} for details, including extension to $\RR^s$ for $s>1$, formulation of two other families of SR wavelets, and the top-down method for point-mass extraction.\\

The consideration of the data function $D({\bf x})$ in (\ref{data fn}) is very natural, since in classical optics, being the Fourier transform of the distribution $h(x)$ in (\ref{pointmass in s-D}), $D({\bf x})$ is the image of $h({\bf y})$, acquired from an ideal thin lens at the focal plane placed at exactly one focal distance behind the lens. On the other hand, the partial differential equation (PDE) that describes the ``photon fluence rate distribution" is the Fokker-Planck PDE, which becomes the isotropic diffusion PDE:
\begin{equation}\label{diffusion}
\begin{cases}
{\displaystyle \frac{\partial}{\partial t} U({\bf x},
t)=c\,\triangledown^2 U({\bf x}, t)},
\;  &{\bf x}\in \RR^s,\,\, t\ge 0;\\
&\\
 {\displaystyle U({\bf x}, 0)= U_0({\bf x})},\, \; &\bf x\in \RR^s,
\end{cases}
\end{equation}
with ``initial value" $U_0({\bf x})$, if the ``absorption" coefficient of photon diffusion of the Fokker-Planck PDE is set to be zero and the isotropic source term is ignored. Here, $c>0$ is the diffusion constant and $\triangledown^2$ denotes the Laplace operator defined by
\begin{equation*}
\triangledown^2  U({\bf x}, t)=\frac{\partial^2}{\partial x_1^2} U({\bf x},t)+ \cdots +\frac{\partial^2}{\partial x_s^2} U({\bf x}, t),
\end{equation*}
where ${\bf x} = (x_1, \cdots, x_s)$. In \cite{C4}, the data function $D({\bf x})$ is modified by scaling the variable $\bf x \in \RR^{s}$ with the parameter ${\bf a} = (a_1,\dots,a_s) \in \RR_{+}^s$, and by applying phase modulation of $D({\bf x})$ with the constant phase  $\kappa >0$ to introduce the ``modified data function" $U_{0}({\bf x};{\bf a})$ in (\ref{diffusion}), used as the ``initial value" of at $t=0$, of the the diffusion PDE with diffusion constant $c>0$. Then the solution $U({\bf x},t;\,{\bf a})$ of this initial-valued PDE is used to isolate the points $\{{\bf x}_{n,k}\}_{k=0}^n$ by choosing a sufficiently large value $t=T>0$. The diffusion process is followed by computing (the complex conjugate of) the inner-product of $U({\bf x},T; \,{\bf a})$ with the wavelet $\Psi_n$ (which is the $s$-fold tensor-product of $\psi_n$ defined in (\ref{Haar wavelets})), after the same phase modulation is applied to $\Psi_n$ and its translation parameter is scaled by $(\frac{1}{a_1},\cdots,\frac{1}{a_s})$ to amplify the coefficients $\{c_{n,k}\}$. This method, introduced and developed in \cite{C4}, is much easier to implement than those developed in \cite{C3}, for solving the inverse problem of recovering the point-mass $\{(c_{n,k}, {\bf x}_{n,k})\}_{k=0}^n$ from the data function $D({\bf x})$. 

\bhag{Epilogue}
 A fair amount of effort in this paper is on the consideration of suitable points that lie on an open arc $\gamma$ in the complex plane as nodes for Lagrange polynomial interpolation. For piece-wise smooth $\gamma$ with corners, in view of the limitation of the exterior conformal map from $|w|>1$ to $\CC^{*}\backslash \gamma$, when extended to $|w|\ge 1$, it is noted that the Fej\'er points on $\gamma$ must be adjusted in order to arrive at the proper growth order for the corresponding Lebesgue constants. An obvious problem for future study is to seek replacement of the Fej\'er points to avoid the use of conformal maps. Besides, when extended to higher Euclidean dimensions $\RR^s$ for $s\ge 3$, a suitable replacement of conformal mapping is hardly available. In this direction, as early as the 1960's, Professor Korevaar already proposed the consideration of Fekete points. However, since the problem of minimizing the total energy $\E({\bf X}_n; D)$ in (\ref {energy2}) is very difficult to attack in the general setting, it seems to be more natural (and more easily approachable) to depend on ``natural phenomena" as opposed to mathematical formulation of ``energy" or ``total energy". In this regard, the ``super-resolution inverse problem" of point-mass recovery from unknown distributions is discussed at the end of Sub-section 5.3. Moreover, it is important to point out that applications of the super-resolution point-mass recovery approach has tremendous real-world and real-life applications. The following list is only a few samples of areas that can be significantly enhanced by incorporating the powerful mathematical tools of super-resolution.
 \vskip.1in
 (1) Fluorescence microscopy: for capturing shapes and colors of live cells  with application to early identification of consequential cancer and pre-canerous change at the earliest stage.
 
 (2) Single-molecule microscopy: for digital detection of bio-molecules for bio-medical applications, including genetic sequencing with single molecules inside a cell, gene editing, and observing coronavirus invasion.
 
 (3) Observational astronomy:  particularly for enhancement of imagery captured by the James Webb space telescope by amplifying the infrared light intensities of the dim imagery of celestial bodies.
 
 (4) Brain surgery: for avoiding accidental damage of nerves by showing their precise location and amplifying their imagery.
 
 (5) Isotope separation: for separating uranium isotope to prepare enriched uranium for use as nuclear fuel, and separating hydrogen isotopes to prepare heavy water for use as moderator in nuclear reactors.
 
 The digital images as discussed above are clusters of point-masses, and can be visualized in any color by using the RGB Bayer format for the point-mass clusters.

\vskip .2in

{\bf{Acknowledgment.}} The first author is indebted to Professor J. Korevaar for the guidance and support during the years of his graduate studies and being a role model over the duration of his academic career. Both authors are grateful to Luyan Wang for his generosity in helping us in the proof of Lemma 2.1.


\begin{thebibliography}{99}
\def\vol#1{{\bf #1}}
\def\sandarbh#1 #2,, #3,, #4..{\bibitem{#1} {\textnormal #2}\ {\it #3}\
#4}
\def\ieee{IEEE International Conference on Neural Networks}
\def\rom#1{\uppercase\expandafter{\romannumeral#1}}
\def\pustak#1{{\rm ``#1'',\/}}

\sandarbh{ABF} E. Abakumov, A. Borichev and K. Fedorovskiy,,, Chui’s conjecture in Bergman spaces,,, Math. Ann. \vol {379} (2021), 1507–1532...

\sandarbh{Belyi1977} V.I. Belyi,,, Conformal mappings and approximation of analytic functions in domains with quasi-conformal boundary,,, Math. USSR-Sb. \vol{31} (1977), 289-317...

\sandarbh{BER1931} S. Bernstein,,, Sur la limitation des valuers d'une polynome,,, Bull. Acad. Sci. USSR, \vol{8} (1931), 1025- 1050...

\sandarbh{B} C. de Boor,,, \pustak{A Practical Guide to Splines},, Applied Mathematical Sciences, \vol {27}, Springer-Verlag (ISBN: 0-387-90356-9),  New York, 1978...

\sandarbh{LB3} L. Bos,,, On Fekete points for a real simplex,,, To appear in this special issue...

\sandarbh{C} C.K. Chui,,, Bounded approximation by polynomials whose zeros lie on a circle,,, Trans. Amer. Math. Soc. \vol{138}, (1969), 171–182...

\sandarbh{C1} C.K. Chui,,, Bounded approximation by polynomials with restricted zeros,,, Bull. Amer. Math. Soc., \vol{73} (1967), 967-972...

\sandarbh{C2} C.K. Chui,,, \pustak { Introduction to Wavelets},, Academic Press--Elsevier (1992), 264 pages...

\sandarbh{C3} C.K. Chui,,, Super-resolution wavelets for recovery of arbitrarily close point-masses with arbitrarily small coefficients,,, Appl. Comput. Harmon. Anal., \vol{61} (2022), 202--253...

\sandarbh{C4} C.K. Chui,,, A diffusion + wavelet-window method for recovery of super-resolution point-masses with application to single-molecule microscopy and beyond,,,  Appl. Comput. Harmon. Anal.,\vol{63} (2023), 1--19...

\sandarbh{C5} C.K. Chui,,, A lower bound of electrostatic fields due to unit point masses,,, Amer. Math. Monthly \vol{78} (1971), 779–780...

\sandarbh{C6} C.K. Chui,,, On approximation in the Bers spaces,,, Proc. Amer. Math. Soc. \vol{40} (1973), 438–442... 

\sandarbh{CS} C.K. Chui and X.C. Shen,,, Order of approximation by electrostatic fields due to electrons,,, Constr. Approx. \vol{1} (1985), 121–135...

\sandarbh{CZ1999} C.K. Chui and L. Zhong,,, Polynomial interpolation and Marcinkiewicz-Zygmund inequalities on the unit circle,,, J. Math. Anal. Appl., \vol{233} (1999),  387--405...

\sandarbh{CZ2021} C.K. Chui and L. Zhong,,, On Marcinkiewicz-Zygmund inequalities and $A_p$-weights for $L$-shape arcs,,, J. Geom. Anal. (2021). https://doi.org/10.1007/s12220-021-00669-2...

\sandarbh{CUR1935} J.H. Curtiss,,,  Interpolation in regularly distributed points,,, Trans. Amer. Math. Soc., \vol{38} (1935), 458–473...

\sandarbh{D} D.L. Donoho,,, Super-resolution via sparsity constraints,,, SIAM J. Math. Anal. \vol{235} (1992), 1309--1331...

\sandarbh{D1977} V.K. Dzjadyk,,, \pustak{Introduction to the Theory of Uniform Approximation of Functions by Polynomials},, Nauk, Moscow, 1977...

\sandarbh{JE} J.M. Elkins,,, Approximation by polynomials with restricted zeros,,, J. Math. Anal. Appl., \vol{25} (1969), 321-336...

\sandarbh{ERD1958} P. Erd\"os,,, Problems and results on the theory of interpolation. I,,, Acta Mathematica Academiae Scientiarum Hungarica, \vol{9} (1958), 381–388...

\sandarbh{F1914} G. Faber,,, \"Uber die interpolatorische Darstellung Stetiger Funktionen,,, Jahresber. Deutsch. Math. Verein., \vol{23} (1914),191-200...

\sandarbh{FE1918} L. Fej\'er,,,  Interpolation und konforme Abbildung,,, Nachr. Ges. Wiss. G\"ottingen. Math.-Phys. Kl. (1918), 319-331... 

\sandarbh{Gaier1987} D. Gaier,,, \pustak{Lectures on Complex Approximation},, Springer, New York, 1987...

\sandarbh{Ib2016} B.A. Ibrahimoglu,,, Lebesgue functions and Lebesgue constants in polynomial interpolation,,, Journal of Inequalities and Applications, (2016), 2016:93 DOI 10.1186/s13660-016-1030-3...

\sandarbh{J} D. Jackson,,, \pustak{The theory of Approximation},, AMS Colloquium Publications, 1930...

\sandarbh{K1} J. Korevaar,,, Asymptotically neutral distributions of electrons and polynomial approximation,,,  The Annals of Mathematics, Second Series, \vol{80} (1964), 403-410...

\sandarbh{K2} J. Korevaar,,, Limits of polynomials with restricted zeros in  \pustak {Studies in Mathematical Analysis and Related Topics.} (Essays in honor of G. Polya),,, Stanford Univ. Press, 1962...

\sandarbh{K3} J. Korevaar,,, Chebyshev-type quadratures: use of complex analysis and potential theory,,, in \pustak{Complex Potential Theory}   P. Gauthier and G. Sabidussi, eds., Dordrecht: Kluwer (1994), 325--364...

\sandarbh{KC} J. Korevaar and C.K. Chui,,, Potentials of families of unit masses on disjoint Jordan curves,,, in \pustak{Abstract Spaces and Approximation}   P.L. Butzer and B-S Nagy eds., Birkhauser Verlag (1969), 338--350...

\sandarbh{KM1} J. Korevaar and J.L.H. Meyers,,, Spherical Faraday cage for the case of equal charges and Chebyshev-type quadrature on the sphere,,, Integral Transforms and Special Functions, \vol{1} (1993), 105--117...

\sandarbh{KM2} J. Korevaar and M.A. Monterie,,, Fekete Potential and Polynomials for Continua,,, J. Approx. Theory, \vol(109) (2001), 110--125...

\sandarbh{KOV1966} T. K\"ovari and C. Pommerenke,,, On Faber polynomials and Faber expansions,,, Math.  Z. \vol{99} (1967), 193--206...

\sandarbh{KS} A.B.J. Kuijlaars and E.B. Saff,,, Asymptotics for minimal discrete energy on the sphere,,, Trans. Amer. Math. Soc. \vol{350} (1998), 523--538...


\sandarbh{MZ} J. Marcinkiewicz and A. Zygmund,,, Mean values of trigonometrical polynomials,,, Fund. Math. \vol{28} (1937),131-166...

\sandarbh {GM} G.R. MacLane,,, Polynomials with zeros on a rectifiable Jordan curve,,, Duke Math. J. \vol{16} (1949), 461-477... 


\sandarbh{DN} D.J. Newman,,, A lower bound for an area integral,,, Amer. Math. Monthly, \vol {79} (1972),  1015–1016...


\sandarbh{P1} D.T. Piele,,, Asymptotically neutral families in $E^3$,,, SIAM J. Math Anal. \vol{4}(1973), 260-268...

\sandarbh{P2} D.T. Piele,,, An Approximation of harmonic functions in $E^3$ by potentials of unit charges,,, SIAM J. Math. Anal. \vol {5} (1974), 563--568...

\sandarbh{R1969} T.J. Rivlin,,, \pustak{An Introduction to the Approximation of Functions},, Blaisdell Publ.Comp. Waltham, Massachusetts, 1969 ...

\sandarbh{RUN1901} C. Runge,,, \"Uber empirische Funktionen und die interpolation zwischen \"aquidistanten Ordinaten,,, Zeitschrift f\"ur Mathematik und Physik, \vol{46} (1901),  224–243...

\sandarbh{SK} E.B. Saff and A.B.J. Kuijlaars,,, Distributing many points on a sphere,,, The Mathematical Intelligence,\vol{19} (1997), 5--11...

\sandarbh{T1} M.D. Thompson,,, Approximation by polynomials whose zeros lie on a curve,,, Duke Muth. J. \vol{31} (1964), 255-265... 

\sandarbh{T2} M.D. Thompson,,, Approximation of bounded analytic functions on the disc,,, Nieuw Arch. Wisk. \vol{159} (1967), 49–54...

\sandarbh{VER1990} P. V\'ertesi,,, Optimal Lebesgue constant for Lagrange interpolation,,, SIAM J. Numer. Anal., \vol{27} (1990), 1322–1331...


\sandarbh{ZZ} L. Zhong and L.Y. Zhu,,, The Marcinkiewicz-Zygmund inequality on a smooth simple arc,,, J. Approx. Theory, \vol{83} (1995), 65--83...

\sandarbh{Z1995} L. Zhong,,, Lagrange interpolation polynomials in $E^p(D)$ with $1 < p < +\infty$,,, Acta Math. Hungar., \vol{69} (1995), 55--66...


\end{thebibliography}
\end {document}